\documentclass[11pt,a4paper]{article}

\providecommand{\keywords}[1]{\par\medskip\noindent\textbf{\textit{Keywords:}} #1}

\usepackage[english]{babel}
\usepackage[T1]{fontenc}
\usepackage[utf8]{inputenc}

\usepackage[
  a4paper,
  top=2.6cm,
  bottom=2.6cm,
  left=2.6cm,
  right=2.6cm
]{geometry}

\usepackage{amsmath,amssymb,amsfonts,mathtools}
\usepackage{graphicx}
\usepackage{tikz}
\usetikzlibrary{positioning,calc}
\usepackage{placeins}
\usepackage{lmodern}
\usepackage{amsmath,amssymb,amsfonts,amsthm}
\usepackage{mathtools}
\usepackage{booktabs}
\usepackage{array}
\usepackage{longtable}
\usepackage{geometry}
\usepackage{xcolor}
\usepackage{lscape,threeparttable}
\usepackage{microtype}
\usepackage{fancyhdr}
\usepackage{titlesec}
\usepackage{enumitem}
\usepackage{algorithm}
\usepackage{algpseudocode}
\usepackage{caption}
\usepackage{subcaption}
\usepackage{url}
\usepackage{pgfplots}
\usepackage{multirow}
\pgfplotsset{compat=newest}
\usepackage{microtype}

\usepackage[colorlinks=true,linkcolor=blue,citecolor=blue,urlcolor=blue]{hyperref}
\usepackage[nameinlink,noabbrev]{cleveref}

\newcommand{\Sset}{\mathcal{S}}
\newcommand{\AinvSelInv}{A^{-1}_{\rm SelInv}}
\newcommand{\AinvInvA}{A^{-1}_{\rm NInv}}
\newcommand{\AinvNInv}{A^{-1}_{\rm NInv}}
\newcommand{\AinvMix}{A^{-1}_{\rm Mix}}
\newcommand{\AinvMixSPAI}{A^{-1}_{\mathrm{Mix\text{-}SPAI}}
}
\newcommand{\AinvCurrent}{A^{\rm inv}}
\newcommand{\AinvCorr}{A_k^{\rm inv}}
\newcommand{\Eop}{E_{\rm inv}}
\newcommand{\invtol}{\delta}
\newcommand{\Lfac}{L}
\newcommand{\Dfac}{D}

\newcommand{\Ainvsel}{\AinvSelInv}
\newcommand{\ANInvpp}{\AinvCurrent}
\newcommand{\Ainvhyb}{\AinvCorr}

\title{Scalable Approximate Selected Inversion Based on Single-Level Incomplete $LDL^T$ Factorization and Spectral Corrections for Large Sparse Systems}

\author{
Tahamina Akter\thanks{Institute for Numerical Analysis, TU Braunschweig, Braunschweig, Germany (\texttt{tahamina.akter@tu-braunschweig.de})}
\and Matthias Bollh\"{o}fer\thanks{Institute for Numerical Analysis, TU Braunschweig, Braunschweig, Germany (\texttt{m.bollhoefer@tu-bs.de})}
\and Olaf Schenk\thanks{Institute of Computing, Faculty of Informatics, Universit\`a della Svizzera italiana, Switzerland (\texttt{olaf.schenk@usi.ch})}
}

\date{}

\begin{document}
\maketitle

\begin{abstract}
This article introduces four parallel numerical techniques for computing entries of the inverse of large sparse symmetric systems, all grounded in incomplete $LDL^T$ (ILDL) factorizations: (1) the selected inversion method (SelInv), which applies the $LDL^T$ factorization to recover entries of the matrix inverse within the sparsity pattern of the computed factors; (2) an approximate inversion method based on a truncated Neumann series expansion applied to the inverse of the L factor (NInv), providing an alternative at the cost of reduced accuracy; (3) a Mix approximation that merges the best of both SelInv and NInv; and (4) Mix-SPAI, which applies sparse approximate inverse (SPAI) refinement on the output of the Mix method to improve entry-level accuracy. To further improve accuracy while maintaining a stable sparsity pattern, we additionally employ a low-rank correction based on eigenvector updates, providing an alternative to tightening the drop tolerance. We report the performance of the proposed numerical techniques on a comprehensive collection of sparse matrices from scientific and industrial applications.
\end{abstract}


\section{Introduction} \label{sec:introduction} 
\begin{figure}[!t] 
\centering \begin{minipage}{0.22\textwidth} 
\centering \includegraphics[width=\linewidth]{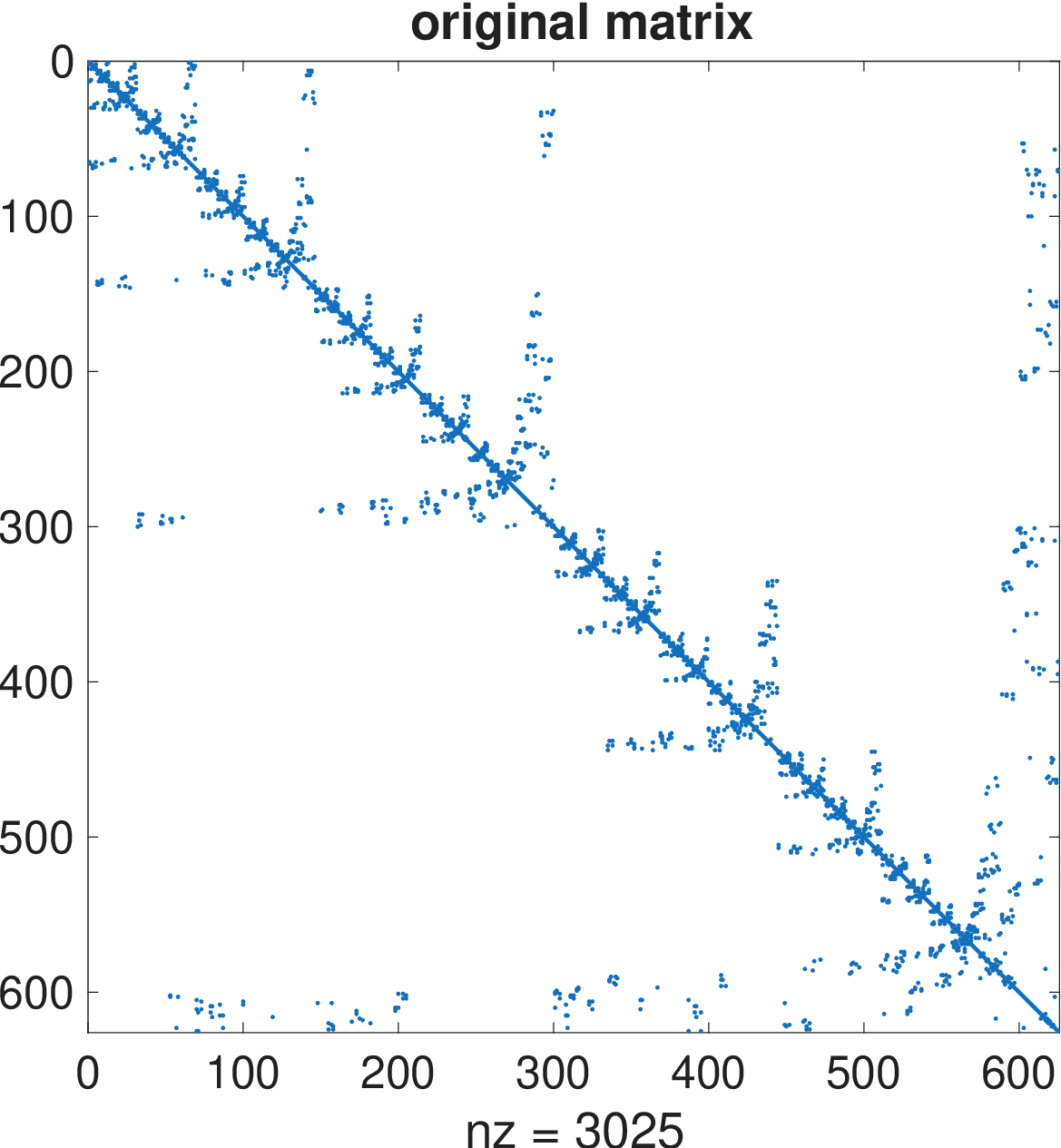} \end{minipage} 
\hspace{0.03\textwidth} 
\begin{minipage}{0.22\textwidth}
\centering \includegraphics[width=\linewidth]{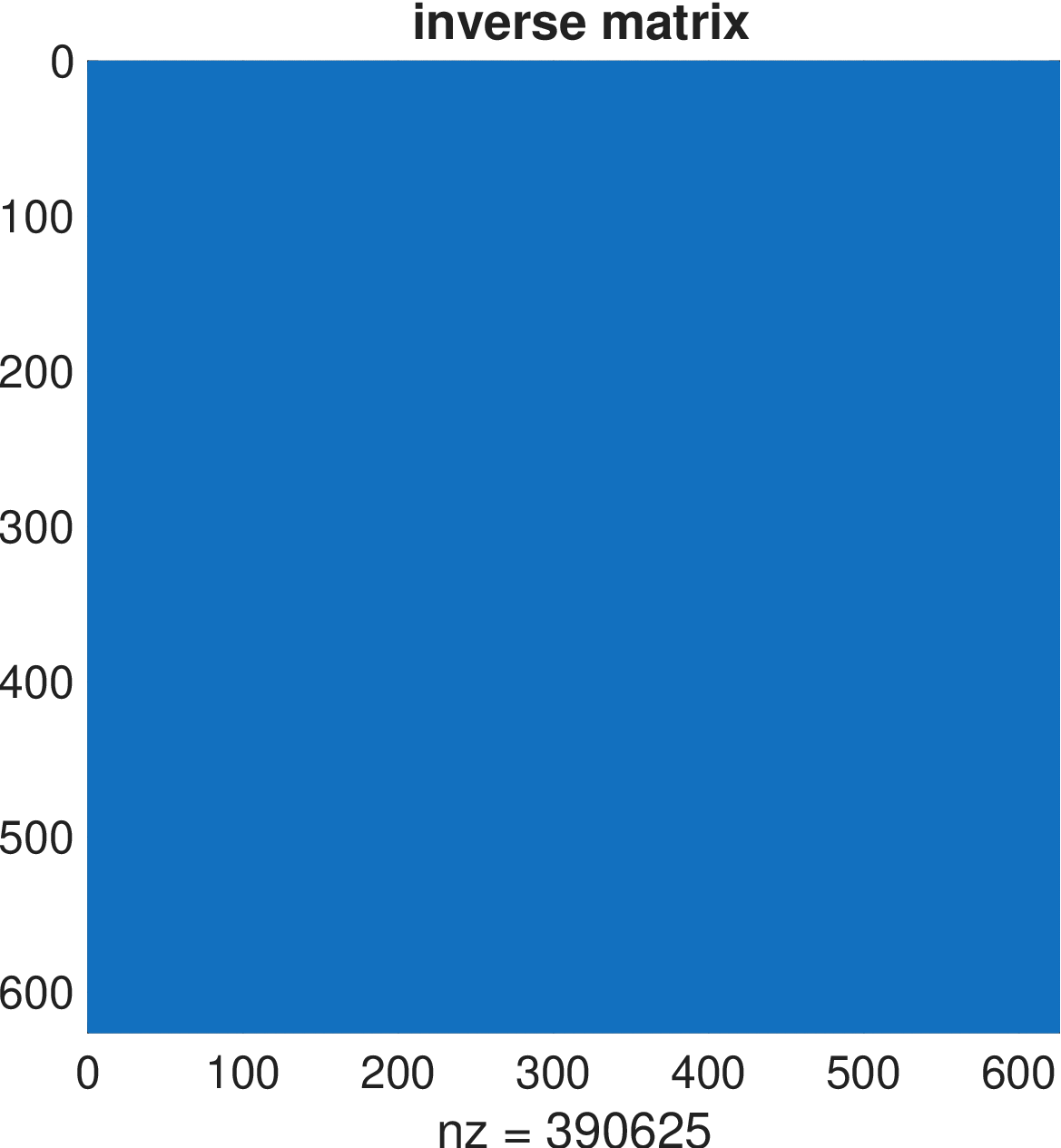} \end{minipage} \hspace{0.03\textwidth} 
\begin{minipage}{0.22\textwidth} \centering \includegraphics[width=\linewidth]{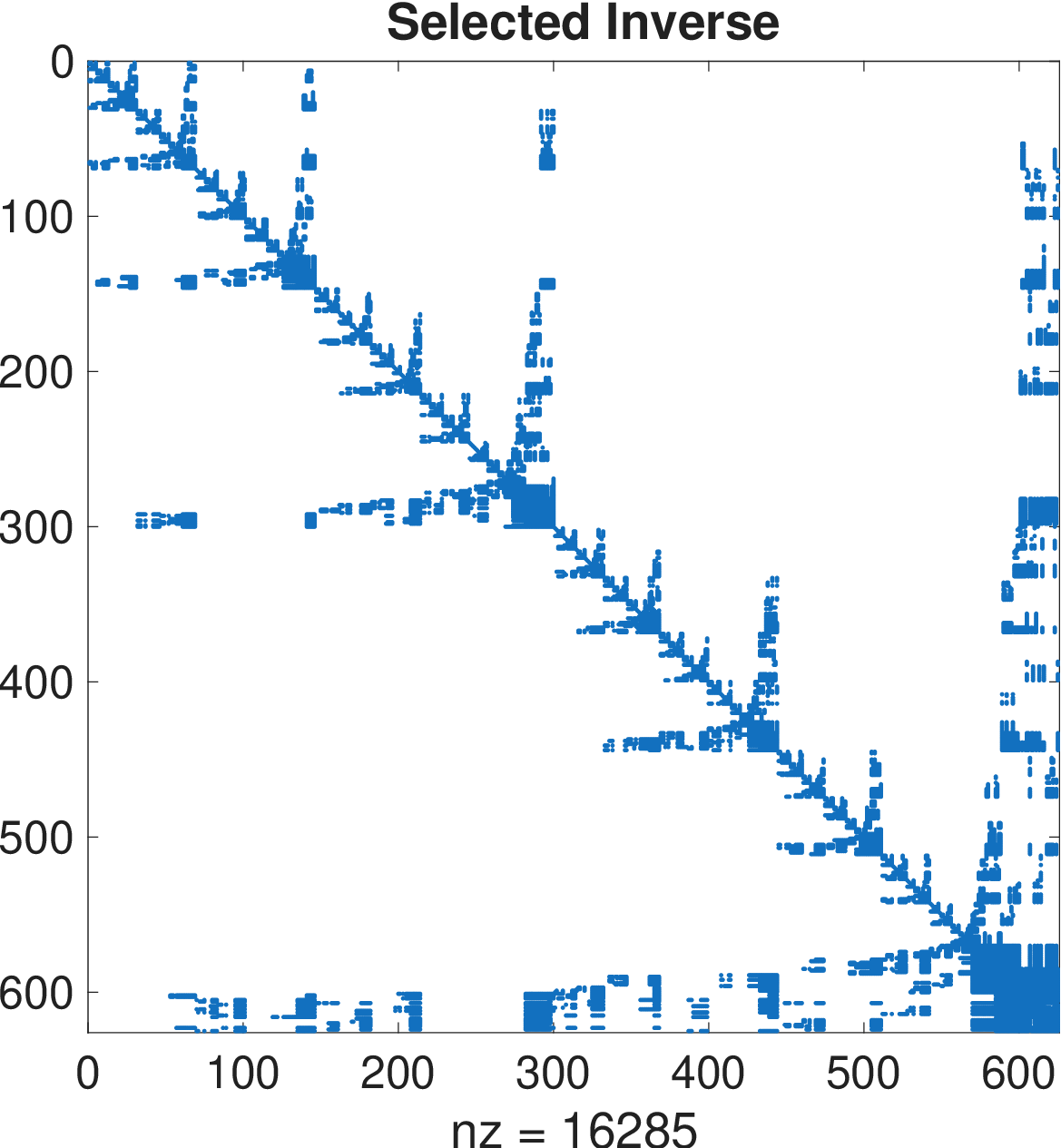}
\end{minipage} 
\caption{Illustration of selected inversion: a sparse matrix $A$ (left), its dense inverse $A^{-1}$ (middle), and a structured subset of the selected inverse entries (right).} \label{fig:matrix_inverse_selinv}
\end{figure}

Matrix inversion, especially for large sparse matrices, appears in many scientific applications and requires high-performance numerical methods \cite{lin_pselinv_2018,8551711,10.1145/3797905.3807841}. Examples include electronic-structure calculations, Green's function evaluation, uncertainty quantification, and inverse covariance matrix estimation. For highly ill-conditioned matrices, forming parts of the inverse remains difficult. In many applications, however, the complete inverse is not required. Instead, only selected entries of $A^{-1}$ are needed. Often, even an approximate inverse is enough. This task is commonly referred to as \emph{selected inversion}. The basic idea is sketched in Figure \ref{fig:matrix_inverse_selinv}: starting from a sparse matrix $A$, one avoids forming the complete inverse and computes only the required entries. Such problems occur in sensitivity analysis, electronic-structure calculations \cite{lin_fast_2011,lin_selinv_2011,rubensson_localized_2021}, quantum transport problems \cite{maillou_serinv_2025} and, discretizations of partial differential equations~\cite{xia_fast_2015}. 

Given  $A\in\mathbb{R}^{n\times n}$, where $A$ is a symmetric, large and sparse matrix, our goal is to compute only entries $(A^{-1})_{ij}$ for an index set $\Sset$ induced by the factorization pattern. To simplify the discussion, suppose that up to a permutation matrix $P\in\mathbb{R}^{n\times n}$, we are given an approximate factorization $P^TAP\approx LDL^T$. Here, we assume that $L$ is block unit lower triangular and $D$ is block diagonal. For ease of presentation, we drop $P$ in the sequel. Given the factors $L,D$, selected inversion uses the pattern
\[ \Sset = \{(i,j):(L+D+L^T)_{ij}\neq 0\}. \]
Thus, the sparsity pattern of the factorization determines which entries of the matrix inverse are computed and stored. The idea of computing selected entries of a matrix inverse goes back to Takahashi et al.~\cite{takahashi_formation_1973}. Lin et al. \cite{lin_selinv_2011} later developed an efficient selected-inversion algorithm based on a sparse factorization. Progress made in this field is closely related to advances in sparse direct solvers and their parallel implementations, including PARDISO~\cite{schenk_pardiso_2004}, SuperLU~\cite{demmel1999superlu}, and MUMPS~\cite{amestoy2001mumps}. The parallel efficiency can be improved further using left-looking selected-inversion methods \cite{jacquelin_left-looking_2018}. Several algorithmic variants have been developed. Lin et al. \cite{lin_selinv_2011} introduced a left-looking parallel selected inversion algorithm for distributed-memory machines and symmetric matrices. Related left-looking approaches have been developed for task-parallel shared-memory systems~\cite{jacquelin_left-looking_2018}, and extensions to nonsymmetric matrices have been considered in PEXSI and related selected-inversion frameworks~\cite{lin_pselinv_2018}. Since diagonal entries and traces are central in many applications, selected inversion has also been combined with techniques for computing diagonal and selected off-diagonal entries of discrete Green's functions \cite{lin_multipole_2009,lin_pole-based_2009}. Direct sparse techniques based on Takahashi's idea have also been used in genomic prediction~\cite{verbosio_enhancing_2017}.

Besides selected inversion, numerical methods for approximately computing the inverse of a matrix provide an alternative way to approximate it.  Techniques of this kind are mostly used as parallel preconditioners to accelerate Krylov subspace methods for solving large-scale linear systems or eigenvalue problems. Among these methods are SPAI-type methods \cite{BarBS99,grote_parallel_1997} which allow for massive parallelization. Another approach, called AINV, directly computes inverse triangular factors without an incomplete factorization \cite{BenCT00}. It has also been successfully used for preconditioning \cite{benzi_comparative_2000}. Alternatively, (incomplete) $LDL^T$ factorization can be employed to compute a sparse approximate inverse by approximately inverting $L$ \cite{BolESS19,eftekhari_block-enhanced_2021}. 
Whether to use a selected inverse approach or an approximate inverse approach is certainly problem-dependent. While the selected inverse applied to an $LDL^T$ factorization computes the inverse within the pattern $\Sset$ exactly, it does not yield any approximation outside. In contrast, computing the matrix inverse approximately yields additional information outside $\Sset$, possibly at the expense of accuracy and memory, particularly when the matrix inverse is relatively dense.

This work presents a framework for scalable approximate inversion and selected-inversion for large sparse matrices. The approach uses an incomplete block $LDL^T$ factorization (ILDL) as the base method, although a direct $LDL^T$ factorization could be used partially. We present two basic inverse approximation methods: a selected inverse and an approximate inverse based on a truncated Neumann series. Since both methods have their strengths and weaknesses, we propose a Mix approximation that merges the best of both approaches. Any of these approximations will be combined with an optional spectral correction based on dominant eigenmodes of the residual inverse-error operator. The goal is to balance accuracy, fill, memory use, and computational time. The implementation uses the JANUS\footnote{JANUS: \url{http://bilu.tu-bs.de}} package, whose underlying methods are based on the level-block approximate $LDL^T$ factorization algorithm in \cite{bollhofer_high_2021}. We build all four inverse constructions on top of the block ILDL factorization, which provides a common computational backbone for SelInv, NInv, Mix, and Mix--SPAI, with the SPAI refinement described in Section~\ref{sec:spai-refinement}.

The rest of this paper is organized as follows: in Section 2, we introduce the necessary mathematical background about matrix inversion and our current methods. Before going into numerical details, we briefly discuss our algorithmic idea. In Section 3, we provide numerical results using some matrices from the SuiteSparse matrix collection\footnote{\url{https://sparse.tamu.edu/}}.  Finally, concluding remarks are drawn in Section 4.

\section{Background and current approaches}
\label{sec:approx-selected-inversion}

Many applications use only a subset of the entries of $A^{-1}$ rather than the full inverse.  
This is because in a large, sparse system, the explicit construction of $A^{-1}$ destroys the sparsity, and it is computationally unreasonable. The aim is to retain the structural integrity of selected entries of the matrix inverse while improving its spectral quality.
Therefore, the purpose of this section is to present a framework for selected inversion and approximate inversion for sparse symmetric systems that combines incomplete factorization 
 and low-rank spectral correction.
We will mainly focus on two approaches, one of which constructs a selected inversion method that builds on prior work on selected inversion-based direct factorizations \cite{lin_selinv_2011,takahashi_formation_1973}. The other one will approximately invert the triangular factors
from the given (incomplete) factorization \cite{eftekhari_block-enhanced_2021}. Throughout this section, unless stated otherwise, we consider a large, sparse, nonsingular, and symmetric matrix $A \in \mathbb{R}^{n\times n}$, and in the main theoretical discussion we assume that there exists an approximate block factorization $A\approx LDL^T$ (up to a symmetric permutation, which we will ignore for simplicity).

\begin{figure}[!t]
    \centering
    \begin{subfigure}[t]{0.25\textwidth}
        \centering
        \includegraphics[width=\linewidth]{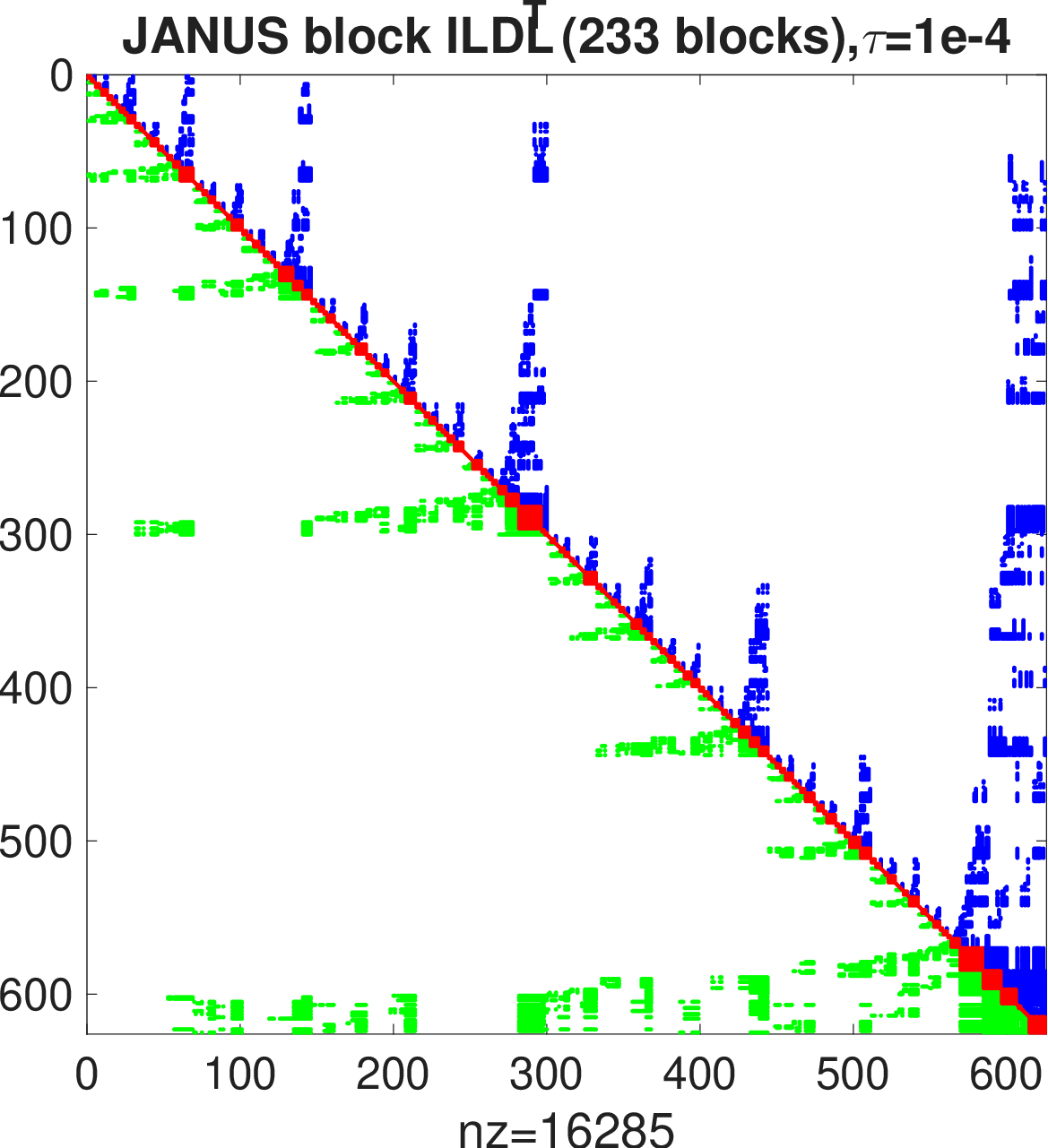}
        \caption{Block-ILDL}
    \end{subfigure}
    \hfill
    \begin{subfigure}[t]{0.25\textwidth}
        \centering
        \includegraphics[width=\linewidth]{fig/selinv.png}
        \caption{Selected inverse}
    \end{subfigure}
    \hfill
    \begin{subfigure}[t]{0.48\textwidth}
      \centering
      \unitlength 1cm
      {
      \begin{picture}(6.75,3.75)
        \put(0.4,0.25){\includegraphics[width=0.8\linewidth]{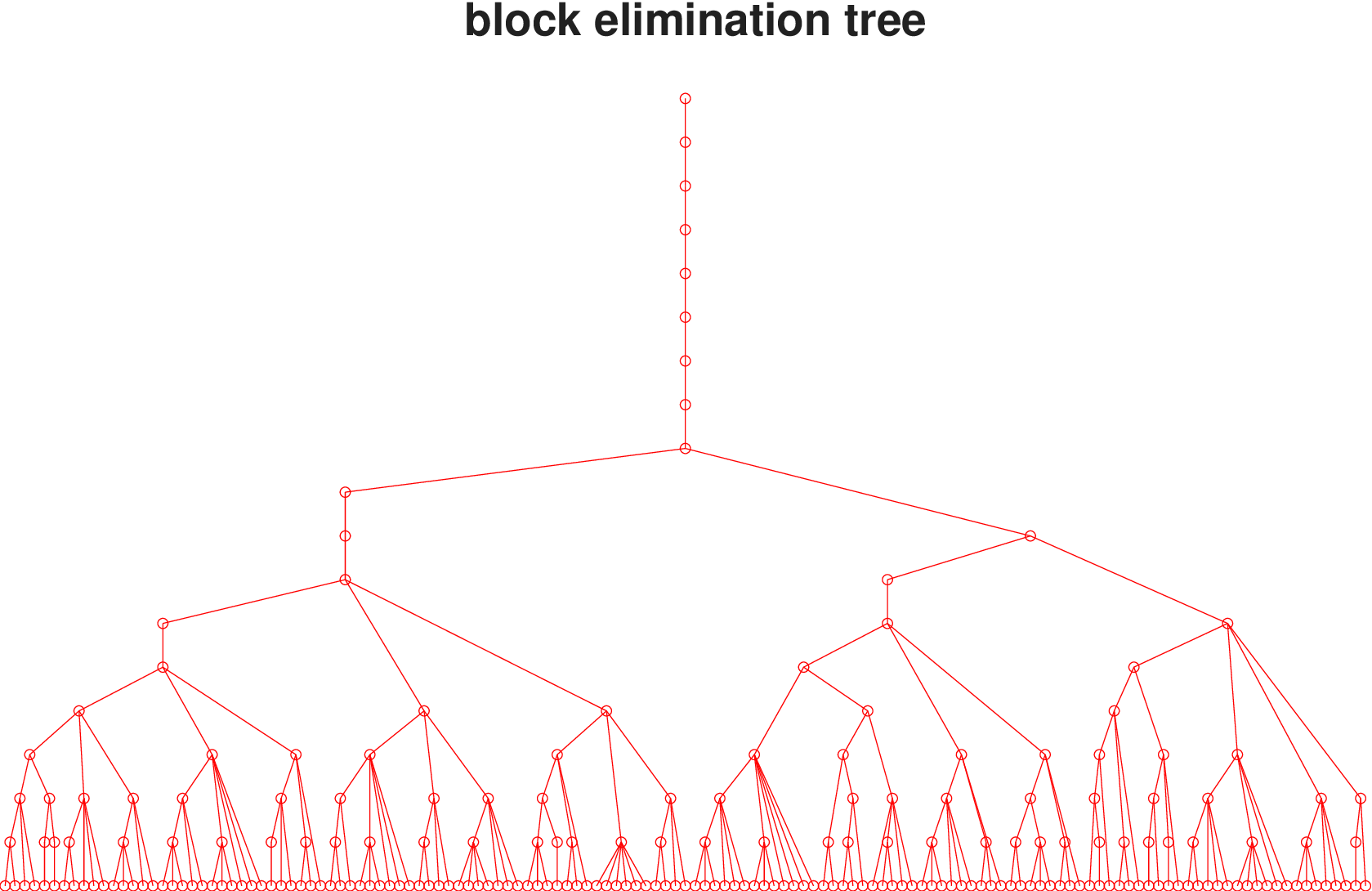}}
        \put(0.25,2.5){$LDL^T \,\uparrow$}
        \put(5.25,2.5){$\downarrow\, \mathrm{SelInv}$}
        \end{picture}
      }
      \caption{Block elimination tree}
    \end{subfigure}
    \caption{Incomplete block factorization (left), approximate selected inverse pattern (middle), and the associated block elimination tree (right).}
    \label{fig:parallel_selinv_overview}
\end{figure}


\subsection{Selected inversion based on incomplete factorization}

At first we consider the selected inversion associated with the block ILDL factorization. A standard way to access selected entries of
$A^{-1}$ is to compute a sparse factorization of $A$. In the symmetric case, we write
\begin{equation}
A\approx LDL^T,
\label{eq:janus_factorization}
\end{equation}
where $L$ is block unit lower triangular and $D$ is block-diagonal with corresponding block structure. The factorization in~\cref{eq:janus_factorization} defines the pattern for the selected inverse. 
In this subsection, we will discuss the selected inverse, denoted by $\AinvSelInv$. In the three succeeding subsections, we present the approximate inverse $\AinvInvA$ based on approximately inverting $L$ and, after that, the Mix approximation $\AinvMix$ and $A^{-1}_{\mathrm{Mix\mbox{-}SPAI}}$. The drop tolerance $\tau$ controls the sparsity and the accuracy of the ILDL. Small values of $\tau$ usually give a more accurate factorization, but also increase fill and runtime. Larger values yield cheaper factors but may lead to a suboptimal approximation.

For the symmetric case, selected inversion algorithms have been developed and implemented in SelInv~\cite{lin_selinv_2011}; distributed-memory parallel variants are described in~\cite{lin_fast_2011,lin_fast_2009}. To recall the general idea of selected inversion, we partition the ILDL in \cref{eq:janus_factorization} as
\begin{equation}
A \approx
\begin{bmatrix}
I & 0\\
L_{21} & I
\end{bmatrix}
\begin{bmatrix}
D_{11} & 0\\
0 & D_{22}
\end{bmatrix}
\begin{bmatrix}
I & L_{21}^T\\
0 & I
\end{bmatrix}.
\label{eq:block_factorization}
\end{equation}
We separate the first leading diagonal block $D_{11}$ of $D$ in the ILDL
factorization from the remaining factorization. $L_{21}$ denotes the
corresponding part of $L$. Since $L$ is block unit lower triangular, the
$(1,1)$ block of $L$ associated with $D_{11}$ must be the identity matrix.

The selected inverse method computes entries of the matrix inverse
only on the set
$\Sset$ induced by the factorization pattern. For the two-by-two block
partitioning in \cref{eq:block_factorization}, inverting
$A$ begins at the bottom right block and proceeds toward to the top
left part.
For a two-by-two step, the selected inverse is obtained via the
following recursion
formula. Suppose that 
$\AinvSelInv$ is partitioned in analogy to \cref{eq:block_factorization}
as
\begin{equation}
\AinvSelInv  =
\begin{bmatrix}
G_{11} & G_{21}^T\\
G_{21} &  G_{22}
\end{bmatrix}
\end{equation}
and suppose furthermore that $G_{22}$ has already been computed.
In this case compute the remaining parts $G_{21}$ and $G_{11}$ via
\begin{equation}
G_{21}\gets -G_{22}L_{21}, 
\label{eq:selinv_recursion1}
\end{equation}
followed by
\begin{equation}
G_{11}\gets D_{11}^{-1}-L_{21}^TG_{21}.
\label{eq:selinv_recursion2}
\end{equation}
To be more precise, the computation of $G_{21}$ is
restricted to the set $\Sset$ and entries
outside this pattern are discarded. This way the computation of $\AinvSelInv$ is recursively traced back to the computation of $G_{22}$ for which a similar relation holds until the final diagonal block of $D$ is reached. Notice that if
$A=LDL^T$ is an exact factorization, then one can show that $\AinvSelInv$
coincides with $A^{-1}$ on $\Sset$.
However, when only an incomplete factorization $A\approx LDL^T$ is available, we could at best expect that $\AinvSelInv$ equals $(LDL^T)^{-1}$ on $\Sset$. Thus, \cref{eq:selinv_recursion1} determines the blocks in which the selected inverse entries are retained according to $\Sset$. Thus, $\AinvSelInv$ preserves the sparsity structure imposed by the incomplete block factorization, while avoiding the construction of the full dense inverse of the system.

In the context of sparse matrix factorization, the inverse information is propagated through elimination trees, along with the structure of the factors, to account for the non-zero pattern of block triangular factors. In the classical selected-inversion algorithm, this propagation proceeds through the elimination tree from parent nodes to child nodes \cite{lin_fast_2011,lin_selinv_2011}. In contrast, the construction of the factorization proceeds in the opposite direction, aggregating information from child nodes to parent nodes, as shown in Figure \ref{fig:parallel_selinv_overview}.


\subsection{Approximate inversion using Neumann series}
\label{subsec:factorized-approx-inversion}
This subsection interprets a method to approximate the inverse of a triangular factor $L$ from the factorization of $A\approx LDL^T$.
We employ an approximate inverse computation method as presented in
\cite{eftekhari_block-enhanced_2021}. We will
refer to this approximate inverse as NInv. The factor 
$L$  can be expressed as  $L = I - E_L$, where  $-E_L$ represents the strict lower triangular part of $L$. 
The exact inverse of  $L$  can be formally expanded as a Neumann series:
\begin{equation}
L^{-1}=I+E_L+E_L^2+E_L^3+\cdots 
\label{eq:neumann}
\end{equation}
Formally, this series is indeed finite since $E_L$ is nilpotent.
In practice, this Neumann-type expansion is truncated and sparsified using
the inverse tolerance parameter $\invtol$. Here, additional powers of $E_L$ are checked versus the computed series with respect to $\invtol$. Eventually we
obtain
\begin{equation}
L^{\rm inv}\approx L^{-1}.
\end{equation}
A larger value of $\invtol$ yields a sparser triangular approximate
inverse. The approximate inverse is then computed via
\begin{equation}
\AinvNInv
\approx
(L^{\rm inv})^TD^{-1}L^{\rm inv},
\label{eq:NInv}
\end{equation}
where again $\invtol$ is used to control the sparsity when forming the
product in \cref{eq:NInv}.
This hopefully gives a sparse
approximation; for more details see \cite{eftekhari_block-enhanced_2021}.
The tolerance is used to control the fill of $L^{\rm inv}$, the approximate
product \cref {eq:NInv} and, therefore, the fill of $\AinvNInv$.

\begin{figure}[!t]
    \centering
    \begin{subfigure}[t]{0.22\textwidth}
        \centering
        \includegraphics[width=0.95\linewidth]{fig/ildl.png}
        \caption{Block-ILDL}
        \label{fig:Block_ildl}
    \end{subfigure}\hfill
    \begin{subfigure}[t]{0.22\textwidth}
        \centering
        \includegraphics[width=0.95\linewidth]{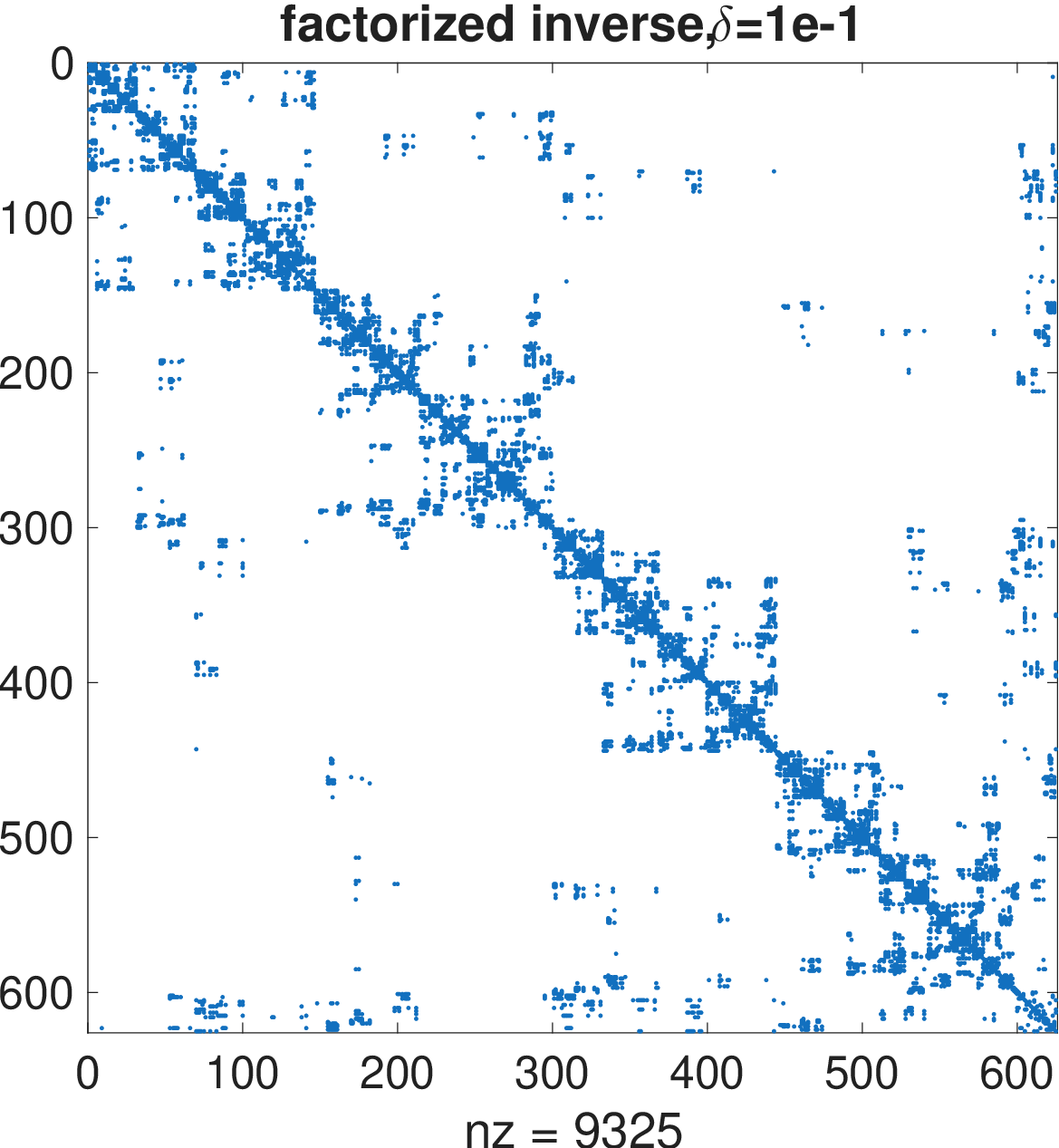}
        \caption{$\invtol=10^{-1}$}
        \label{fig:fspai_ildl_1e-1}
    \end{subfigure}\hfill
    \begin{subfigure}[t]{0.22\textwidth}
        \centering
        \includegraphics[width=0.95\linewidth]{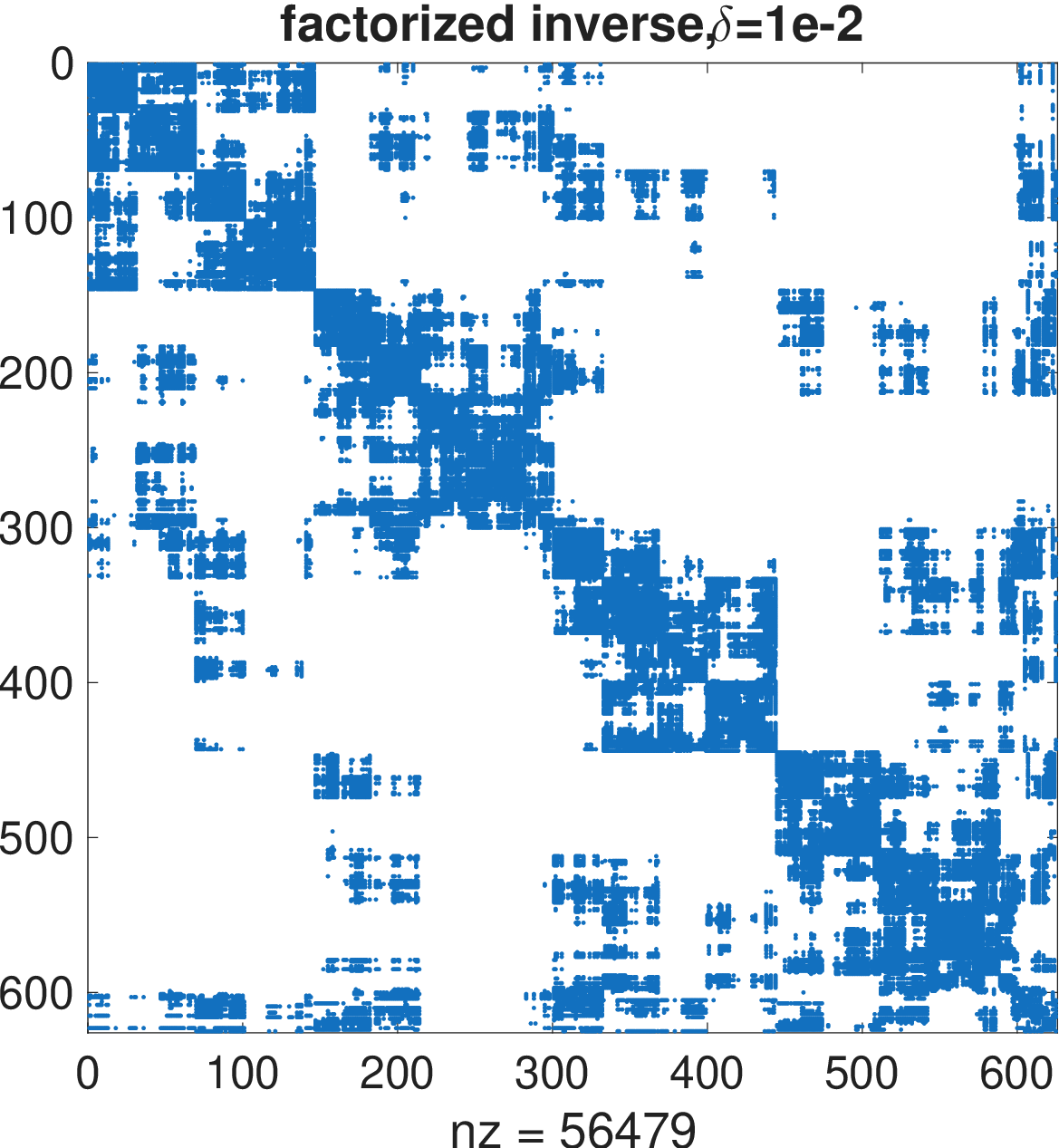}
        \caption{$\invtol=10^{-2}$}
        \label{fig:fspai_ildl_1e-2}
    \end{subfigure}\hfill
    \begin{subfigure}[t]{0.22\textwidth}
        \centering
        \includegraphics[width=0.95\linewidth]{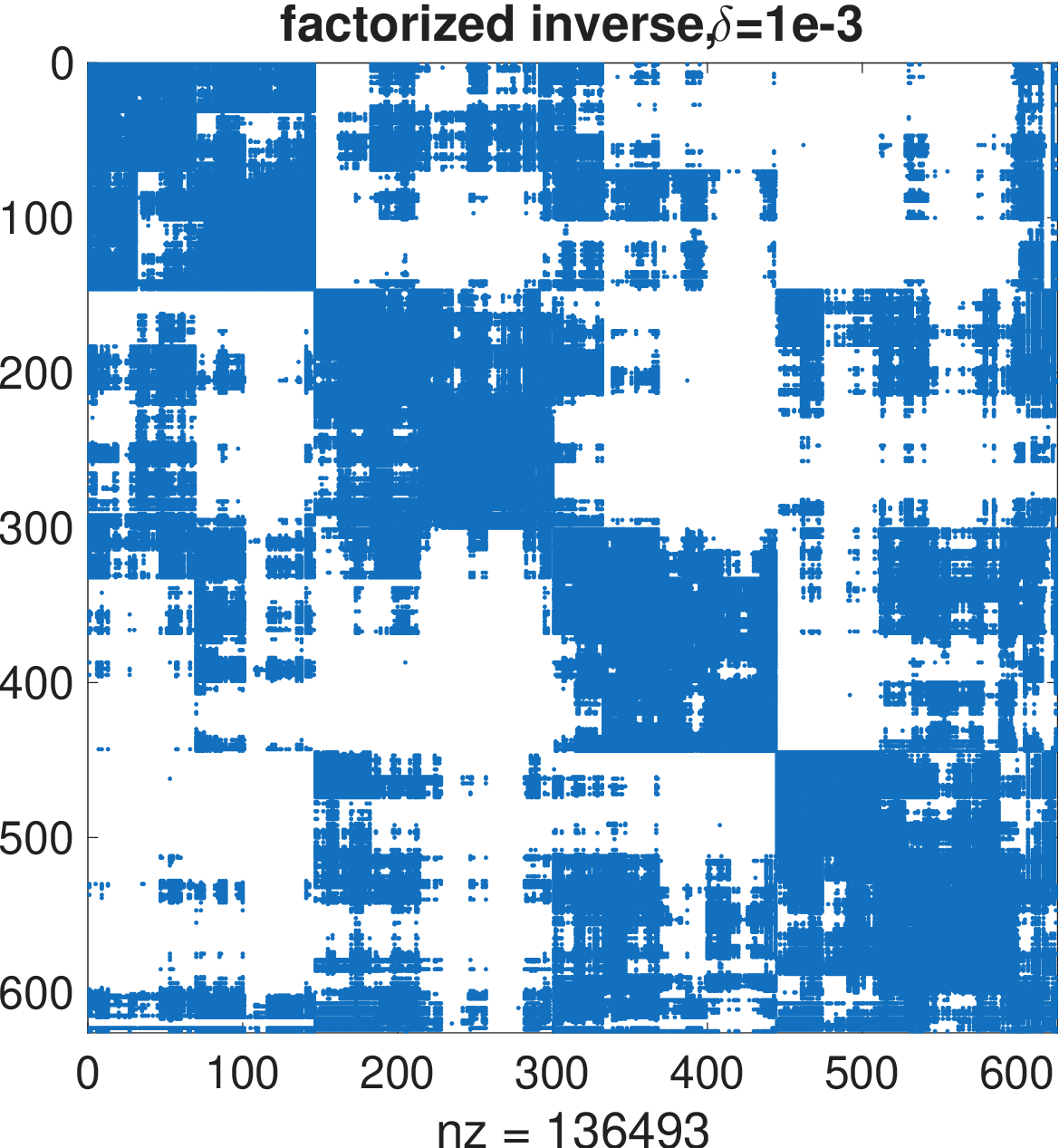}
        \caption{$\invtol=10^{-3}$}
        \label{fig:fspai_ildl_1e-3}
    \end{subfigure}
    \caption{Incomplete block factorization (left) and the corresponding
    factorized approximate inverses for three different inverse tolerances
    $\invtol$ (right three figures). Decreasing $\invtol$ keeps more entries
    in $\AinvNInv$ and leads to a denser approximation.}
    \label{fig:blockildl_fspai_comparison}
\end{figure}

In Figure  \ref{fig:blockildl_fspai_comparison} we demonstrate the impact
of $\invtol$. 
A smaller value of $\invtol$, e.g.,  from $10^{-1}$ to $10^{-3}$ in Figure 
\ref{fig:blockildl_fspai_comparison}, keeps more entries and possesses more
off-diagonal entries. This gives a denser sparsity pattern for
$L^{\rm inv}$ and hence for $\AinvNInv$. A looser tolerance provides a
cheaper approximation, while a tighter tolerance may substantially
increase the fill. Unlike $\AinvSelInv$, which is stored only on the
selected pattern, $\AinvNInv$ demonstrates a sparse approximate inverse where the fill only depends on the size of the entries.

\subsection{Mixed approximate inversion}
\label{subsec:Mix-inversion}

The selected inverse and the approximate inverse contain different information. Each of them has its strength but also its weakness.
The Neumann-series based approximation $\AinvInvA$ gives a global sparse inverse, only
controlled by $\invtol$. In the best case this could lead to an approximation
of $A^{-1}$ on the order of $\invtol$. However, if $A^{-1}$ has many
large entries, then $\AinvNInv$ is likely to be dense.
$\AinvSelInv$ gives selected entries on the pattern of the block ILDL. If
the factorization were exact, then $\AinvSelInv$ would coincide with $A^{-1}$
on the pattern $\Sset$ of the $LDL^T$ factorization. While this kind of inverse could
lead to a high approximation of $A^{-1}$ inside this pattern, no
approximation is computed outside.

This complementary behavior motivates us to merge both approaches to get
the best out of both of them.
In order to balance between these different approaches we will 
introduce a Mix approximation.

Let $\Sset_{\varepsilon}\subseteq\Sset$ be the set of selected entries
such that the entries of $|\AinvSelInv|$ are above a prescribed target tolerance $\varepsilon$. The Mix approximation is defined by
\[
(\AinvMix)_{ij}
=
\begin{cases}
(\AinvSelInv)_{ij}, & (i,j)\in\Sset_{\varepsilon},\\[3pt]
(\AinvNInv)_{ij}, & \text{otherwise}.
\end{cases}
\]
Thus, $\AinvMix$ keeps the sparse structure of $\AinvNInv$, but replaces
selected entries in $\Sset_{\varepsilon}$ by the corresponding values from $\AinvSelInv$.
We note that $\AinvMix$ now depends on three parameters. At first, there
is the drop tolerance $\tau$ for the block ILDL. Second, there is
$\invtol$ which is used for $\AinvNInv$. Third, we have introduced
the parameter $\varepsilon$ for dropping entries of small size in $\AinvSelInv$.
Now in order to avoid that the computation of $\AinvMix$ is at the
price and the fill of the computation of both, $\AinvSelInv$ and $\AinvNInv$,
we use a loose tolerance $\invtol$ when computing $\AinvNInv$, whereas
the size of $\varepsilon$ will be used as target tolerance for the overall
approximate inverse.

In the sequel, all four methods may be improved
further by low-rank corrections as explained in the next section.
When we refer to $\AinvCurrent$, it may denote any of
$\AinvSelInv,\AinvNInv,\AinvMix$ or later also $\AinvMixSPAI$(after being introduced Section \ref{sec:spai-refinement}).


\subsection{SPAI Refinement of Approximate Inverses}
\label{sec:spai-refinement}
The sparse approximate inverse (SPAI) \cite{grote_parallel_1997} refinement is used here as an
additional value-improvement step on a prescribed sparsity pattern. Starting from an already computed approximation
$\AinvCurrent$, we keep its sparsity pattern fixed and improve the nonzero values column by column. For the $j$-th column, SPAI seeks an approximation $m_j$ that reduces
\[
    \|A m_j-e_j\|_2,
\]
subject to the sparsity pattern inherited from the initial
approximation. The corresponding normal equations are solved by diagonally preconditioned conjugate gradients. Since the individual columns are independent, this refinement is naturally parallel. In the fourth method considered in the numerical experiments, denoted Mix--SPAI, the starting approximation is the Mix inverse. First, a relatively inexpensive NInv approximation is constructed using the
loose inverse tolerance (e.g. $\delta=10^{-1}$). Selected entries are then overwritten by the corresponding SelInv values according to the mixed strategy described in Section~\ref{subsec:Mix-inversion}. Finally, SPAI refines the values on the resulting fixed sparsity pattern. Thus,
\[
   A^{-1}_{\mathrm{Mix\mbox{-}SPAI}}
   =
   \operatorname{SPAI}\!\left(A^{-1}_{\mathrm{Mix}}\right),
\]
where SPAI changes the values but does not enlarge the prescribed pattern. This construction differs from NInv, which uses the tighter tolerance $\delta=10^{-2}$ as an independent method. The purpose of the loose NInv stage in Mix--SPAI is to obtain an inexpensive candidate pattern,
while SelInv and SPAI subsequently improve the entries on that pattern.

To maintain symmetry, eventually $\frac12(A^{-1}_{\mathrm{Mix\mbox{-}SPAI}}+(A^{-1}_{\mathrm{Mix\mbox{-}SPAI}})^T)$ is used. 
\subsection{Inversion error and spectral correction}
\label{subsec:error-decomposition}

It is useful to distinguish two sources of error. The first is the factorization error introduced by the block ILDL,
\[
E_{\rm fact}=A-LDL^T.
\]
This error depends on the drop tolerance $\tau$, the ordering, the scaling, and the block structure.

The second error is obtained after computing the approximate inverse $\AinvCurrent$. We define
\[
\Eop=A^{-1}-\AinvCurrent.
\]
Thus, $E_{\rm fact}$ measures the quality of the incomplete
factorization, while $\Eop$ measures the remaining error in the inverse approximation. The spectral correction introduced below acts on $\Eop$.
Suppose we supplement $\AinvCurrent$ with a rank-$k$-correction
$V_k\Lambda_kV_k^T$, where $V_k\in\mathbb{R}^{n\times k}$, $\Lambda_k\in\mathbb{R}^{k\times k}$. We denote the rank-$k$ corrected approximation of $A^{-1}$ by
\begin{equation}
\AinvCorr
:=
\AinvCurrent + V_k\Lambda_kV_k^T.
\label{eq:Ainv-low-rank}
\end{equation}
Then the remaining inverse error is
\[
A^{-1}-\AinvCorr
=
\Eop-V_k\Lambda_kV_k^T.
\]
The main idea for adding a low-rank correction can be explained by the behavior of the incomplete factorization. When using a drop tolerance $\tau$,
the factorization error
$\|E_{\rm fact}\|$ is expected to be controlled approximately by the order of magnitude $\tau\|A\|$.
For symmetric matrices, Weyl's perturbation theorem implies that the absolute perturbation of each eigenvalue is bounded by $\|E_{\rm fact}\|_2$. Consequently, eigenvalues whose magnitude is
large relative to this perturbation scale are typically affected by a smaller relative error, whereas small-magnitude eigenvalues may be considerably more sensitive to the incomplete factorization. However, simply looking at the eigenvalues of small magnitude of $A$ (resp. eigenvalues of large magnitude of $A^{-1}$) would not take into account how well $A^{-1}$
 is approximated by $\AinvCurrent$. The relevant quantity is therefore the remaining inverse-approximation error
\[
\Eop=A^{-1}-\AinvCurrent.
\]
A low-rank correction is particularly useful when the dominant
components of $\Eop$ can be represented accurately by only a small
number $k$ of eigenmodes.

Let the eigenvalues of $\Eop$ be ordered such that
\[
|\lambda_1(\Eop)|
\ge |\lambda_2(\Eop)|
\ge \cdots .
\]
Since $\Eop$ is symmetric, its singular values satisfy
\[
\sigma_i(\Eop)=|\lambda_i(\Eop)|.
\]

Given $\Eop=A^{-1}-\AinvCurrent$, then the optimal rank-$k$ approximation
is obtained from the Eckart--Young--Mirsky theorem theorem for the spectral
norm, since 
\[
\|\left(A^{-1}-\AinvCurrent\right) \ - \ V_k\Lambda_kV_k^T\|_2 = 
\|\Eop  -  V_k\Lambda_kV_k^T\|_2
=
|\lambda_{k+1}(\Eop)|
\]
holds when choosing $\Lambda_k$ as diagonal matrix with the leading $k$ dominant eigenvalues $\lambda_{1}(\Eop),\dots,\lambda_{k}(\Eop)$ and when their associated
orthonormal eigenvectors are chosen as columns of $V_k$.
The same arguments hold for the Frobenius norm, since for this choice
we also get the optimal bound
\[
\|\Eop  -  V_k\Lambda_kV_k^T\|_F^2
=
\sum_{i>k}\lambda_i(\Eop)^2.
\]
For perturbation bound analysis, we refer the reader to \cite{horn_matrix_2012}.

Although the optimal rank-$k$ approximation can be stated explicitly
by the dominant eigenspace $(V_k,\Lambda_k)$ of $\Eop$,
 \[
 \Eop V_k= V_k\Lambda_k,
 \qquad
 V_k^TV_k=I,
 \]
its numerical computation formally
requires the exact $A^{-1}$ which is not available.
Any method
for computing a low-rank correction from $\Eop$ has to replace $A^{-1}$
by some approximation which is more accurate than $\AinvCurrent$.
A natural way to deal with $A^{-1}$ will be to use a matrix-free
method. In that case, for a vector $x$, we apply
\[
\Eop x=A^{-1}x-\AinvCurrent x.
\]
We will substitute the first term $A^{-1}x$ by a
preconditioned Krylov subspace method where the existing
block ILDL factorization is used as a preconditioner.
The Krylov subspace methods used here are the preconditioned CG method \cite{golub_matrix_2013}
in the positive definite case, respectively the preconditioned simplified QMR
method \cite{FreN95} in the indefinite case.
For the second term $\AinvCurrent x$, we can easily apply $\AinvCurrent$ directly, since $\AinvCurrent$ is explicitly available.

Thus, the correction removes the dominant unresolved modes of the inverse error without increasing the sparse fill of the ILDL factor. The overall computational procedure is summarized in Algorithm \ref{alg:main}.

\begin{algorithm}[ht!]
\caption{Meta-algorithm: Block Approximate Inverse with Spectral Correction}
\label{alg:main}
\begin{algorithmic}[1]
  \State Prescribe $\tau,\delta,\varepsilon,k$
  \State \textbf{Block factorization.} Compute a block incomplete factorization
$A\approx LDL^T$ with given tolerance $\tau$.

\State \textbf{Base inverse approximation.} Construct a base inverse
approximation $\AinvCurrent\approx A^{-1}$, where
$\AinvCurrent\in\{\AinvSelInv,\AinvNInv,\AinvMix,\AinvMixSPAI\}$ according to
$\delta$, $\varepsilon$ where needed.

\State \textbf{Residual inverse operator.} Define the matrix-free
inverse-error operator $\Eop=A^{-1}-\AinvCurrent$.

\State \textbf{Spectral correction.} Compute a few dominant eigenpairs
$\Eop V_k\approx V_k\Lambda_k$ and form
$\AinvCurrent_k=\AinvCurrent+V_k\Lambda_kV_k^T$.

\State \textbf{Acceptance test.} Use the corrected diagonal only if it
improves the diagonal error; otherwise keep the base approximation
$\AinvCurrent$.
\end{algorithmic}
\end{algorithm}

Before going into the details of the numerical discussion, we state the components of the measurement.

\subsection{Measured quantities}
Since for large-scale systems we can hardly provide the exact inverse,
we will measure the quality of our approximate inverses by measuring the
diagonal entries of the inverse. The exact values will be computed
in advance by a direct solver. Our approximate inverse will not
make use of this information.
Let $V_k=\left(v_{ij}\right)_{i=1,\dots,n, j=1,\dots,k}\equiv \left(v_1,\dots,v_k\right)$. 
If we denote by
\[
\widehat d=\operatorname{diag}(\AinvCurrent),
\]
then the diagonal entries of $\AinvCorr$ are given by 
\[
\widehat d^{k}
=
\operatorname{diag}(\AinvCorr)
=
\widehat d+
\sum_{j=1}^{k}\lambda_j(v_j\odot v_j),
\]
or, entrywise,
\[
(\widehat d^{k})_i
=
\widehat d_i+
\sum_{j=1}^{k}\lambda_j v_{ij}^2, \; i=1,\dots,n.
\]
The low-rank correction is stored in factored form and requires only
$\mathcal{O}(nk)$ additional storage and work, with $k\ll n$.

This correction is used only adaptively. Tightening the drop tolerance $\tau$ reduces the factorization error, but also increases fill and runtime. So, eigenvector updates are useful only when the dominant eigenvectors of $\Eop$ have a visible effect on the target quantities. In the numerical experiments, we
therefore keep the correction only when it improves the diagonal error.

We also set
\[
d=\operatorname{diag}(A^{-1}).
\]
The relative diagonal error is measured by
\[
\operatorname{err}_1
=
\frac{\|d-\widehat d\|_1}{\|d\|_1}.
\]
When the corrected diagonal is accepted, $\widehat d$ is replaced by
$\widehat d^{k}$ in the same expression.

We also report the factor fill and inverse fill,
\[
\operatorname{fill}_{\rm fact}
=
\frac{\operatorname{nnz}(L+D+L^T)}{\operatorname{nnz}(A)},
\qquad
\operatorname{fill}_{\rm inv}
=
\frac{\operatorname{nnz}(\AinvCurrent)}{\operatorname{nnz}(A)}.
\]
The runtime is split as
\[
T_{\mathrm{total}}
=
T_{\mathrm{fact}}
+
T_{\mathrm{inv}}
+
T_{\mathrm{eig}},
\]
 
where $T_{\rm fact}$ is the factorization time, $T_{\rm inv}$ is the time required to construct the base inverse approximation, and $T_{\rm eig}$ is the time spent in the spectral correction.

\section{Numerical Experiments}
\label{sec:numerical-experiments}

In this section, the proposed approximate inverse framework is investigated with a set of large sparse symmetric test matrices from the SuiteSparse Matrix Collection\footnote{\url{https://sparse.tamu.edu/}}. The experiments are designed to assess four aspects of the framework: the sensitivity with respect to the block ILDL with drop tolerance $\tau$, the accuracy and sparsity of the different base inverse approximations, the numerical effect of the spectral correction, and the distribution of the total runtime among factorization, inverse construction, and eigensolution. In addition to SelInv, NInv, and Mix, we include Mix--SPAI, in which SPAI refinement \cite{grote_parallel_1997} is applied to the fixed
sparsity pattern of the Mix approximation.

\subsection{Experimental and computational setup}
All large-scale numerical experiments have been performed on a Linux-based system using MATLAB interfaces. The node used in the experiments has 1~TB of main memory and 4 Intel(R) Xeon(R) CPU E7-4880 v2 processors at 2.50~GHz, with 15 cores per socket and 60 cores in total. For the block ILDL factorization, we use the method from \cite{bollhofer_high_2021} as implemented in the JANUS\footnote{JANUS: \url{http://bilu.tu-bs.de}} package. Our selected inverse method to compute $\AinvSelInv$ is built on top of it. It is parallelized in analogy to existing methods based on the elimination tree (cf. \cite{lin_selinv_2011}). For the approximate inverse $\AinvNInv$, we follow the implementation described in \cite{eftekhari_block-enhanced_2021}, which is part of the SQUIC\footnote{SQUIC: \url{https://www.gitlab.ci.inf.usi.ch/SQUIC}} package and is also parallelized. The dominant eigenpairs of the residual inverse-error operator are computed using PRIMME\footnote{PRIMME: \url{https://github.com/primme/primme}} as presented in \cite{stathopoulos_primme_2010}. Additional experiments with MATLAB \texttt{eigs} and JaDaMILU\footnote{JaDaMILU: \url{http://metronu.ulb.ac.be/JADAMILU/}} as published in \cite{m_bollhofer_jadamilu_2007} gave qualitatively consistent results, while PRIMME generally required less memory and computational time for the matrices considered
here. Unless stated otherwise, all experiments use
\[
k=n_{\mathrm{eig}}=3,
\qquad
\tau\in\{10^{-3},10^{-4},10^{-5},10^{-6}\},
\]
where $\tau$ denotes the drop tolerance of the block ILDL
factorization. For NInv, the factorized approximate inverse is
constructed with $\delta=10^{-2}$. For Mix, a looser tolerance
$\delta=10^{-1}$ is used for the NInv component before selected entries are replaced by their SelInv values. Mix--SPAI starts from the same Mix construction with $\delta=10^{-1}$ and subsequently applies SPAI refinement on the fixed Mix sparsity pattern.
These parameter choices intentionally favor inexpensive approximate inverses, allowing the selected inverse and the subsequent spectral correction to recover most of the lost accuracy while maintaining a lower computational cost.

All cross-method comparisons in this section are taken from the same selected-suite experiment. We note that our ILDL method is a sequential algorithm. Therefore,
its computation time is sometimes dominating compared to the computation time for the approximate inverse methods for building $\AinvSelInv$ and $\AinvNInv$, which are both parallelized. The eigenvector computation time with PRIMME is also relatively large because not all parts are parallelized. 
In principle, matrix-vector products and vector updates expose fine-grained parallelism \cite{chow_fine-grained_2015}. 
In the present implementation, this parallelism is mainly obtained through the sparse solves and the application of $\AinvCurrent$ \cite{chow_using_2018}. 

\subsection{Properties of Test Matrices}

The matrix characteristics of the representative test matrices are listed in Table \ref{tab:testmatrices}.
Here $n$ denotes the matrix dimension and $\operatorname{nnz}(A)/n$ is the average number of nonzeros per row. We also report the degree of diagonal dominance,
\[
\operatorname{DoDD}(A)
=
\min_i
\frac{|a_{ii}|}{\sum_{j\neq i}|a_{ij}|}
\]
This quantity gives a simple structural indication of how strongly the diagonal dominates the off-diagonal entries. Values from Table \ref{tab:testmatrices} indicate that at least one row is not strictly diagonally dominant. In our experiments, covariance-type matrices with stronger diagonal dominance are usually well approximated by the base inverse approximation, while structural matrices such as \texttt{crankseg\_1} and \texttt{crankseg\_2} may require spectral correction to improve the diagonal accuracy. Related perturbation results for diagonally dominant matrices can be found in \cite{dailey_relative_2014}. The microbenchmark matrices are substantially denser than the representative SuiteSparse matrices. In particular, the matrices 
\texttt{NEW\_A\_microbench\_n12\_b1024}, and\linebreak \texttt{NEW\_A\_microbench\_n24\_b1024} have very large values of \(\operatorname{nnz}(A)/n\), which makes them useful for evaluating the scalability of the parallel quadratic selected inversion \cite{10.1145/3797905.3807841} as well as for our application. The two microbenchmark matrices are strongly diagonally dominant according to the reported DoDD values
; see Table \ref{tab:microbench_matrices}.

\begin{table}[!t]
\centering
\caption{Sample matrix information used in the numerical experiments, sorted by matrix dimension \(n\).}
\label{tab:testmatrices}
\begin{tabular}{lrrll}
\toprule
Matrix & $n$ & $\operatorname{nnz}(A)/n$ & Definiteness & DoDD \\
\midrule
bodyy6          & 19,366  & 7.0   & Indefinite & $1.0$ \\
wathen120       & 36,441  & 15.5  & SPD        & $1.7{\times}10^{-1}$ \\
bcsstk39        & 46,772  & 44.7  & SPD        & $8.1{\times}10^{-2}$ \\
bcsstm39        & 46,772  & 1.0   & Indefinite & $\infty$ \\
crankseg\_1     & 52,804  & 201.0 & SPD        & $2.3{\times}10^{-2}$ \\
crankseg\_2     & 63,838  & 221.6 & SPD        & $2.2{\times}10^{-2}$ \\
oilpan          & 73,752  & 48.8  & Indefinite & $1.5{\times}10^{-1}$ \\
apache1         & 80,800  & 6.7   & SPD        & $1.0$ \\
2cubes\_sphere  & 101,492 & 16.2  & SPD        & $3.3{\times}10^{-1}$ \\
bmwcra\_1       & 148,770 & 71.5  & Indefinite & $1.4{\times}10^{-1}$ \\
pwtk            & 217,918 & 53.4  & SPD        & $1.8{\times}10^{-5}$ \\
\bottomrule
\end{tabular}
\end{table}

\begin{table}[!t]
\centering
\caption{Microbenchmark matrix information taken from \cite{10.1145/3797905.3807841}.}
\label{tab:microbench_matrices}
\begin{tabular}{lrrll}
\toprule
Matrix & $n$ & $\operatorname{nnz}(A)/n$ & Definiteness & DoDD \\
\midrule
NEW\_A\_microbench\_n12\_b1024   & 12,288  & 1963.6 & SPD        & $1.8$ \\
NEW\_A\_microbench\_n24\_b1024   & 24,576  & 2006.3 & SPD        & $1.8$ \\
\bottomrule
\end{tabular}
\end{table}

\subsection{Effect of the Drop Tolerance}

The drop tolerance $\tau$ is the primary parameter governing the block ILDL factorization, the core part of the computation. Smaller values of $\tau$ retain more entries in the factors and usually improve the inverse approximation, but they also increase memory consumption and runtime. Larger values of $\tau$ give sparser factors, but may lead to a larger inverse approximation error,
\[
\Eop=A^{-1}-\AinvCurrent.
\]
Thus, the choice of $\tau$ directly affects both the fill-in and the inverse accuracy.



\begin{figure}[!t]
    \centering
    \includegraphics[width=\textwidth]{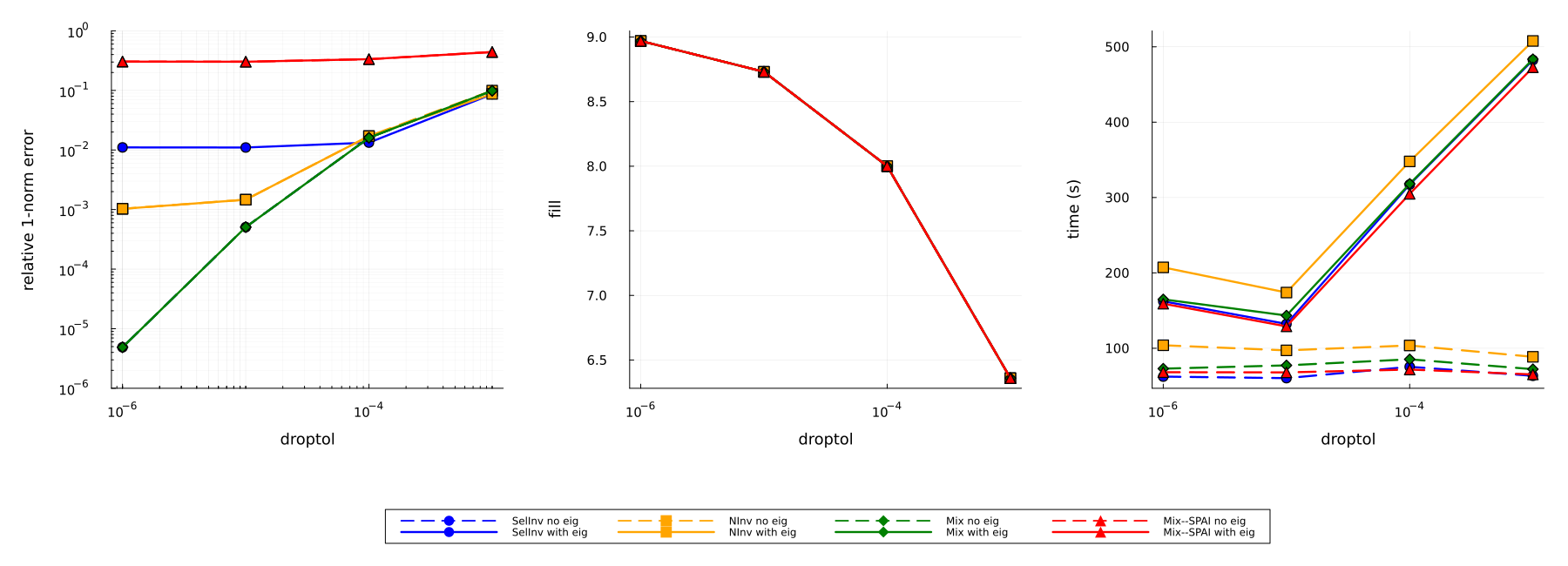}
    \caption{Detailed comparison of SelInv, NInv, Mix, and Mix--SPAI
    for the matrix \texttt{crankseg\_2}. The panels show the relative
    diagonal error with and without spectral correction, fill from the ILDL factor, and
    the runtime components as functions of the drop tolerance.}
    \label{fig:crankseg2_detailed}
\end{figure}

\begin{figure}[!t]
    \centering
    \includegraphics[width=\textwidth]{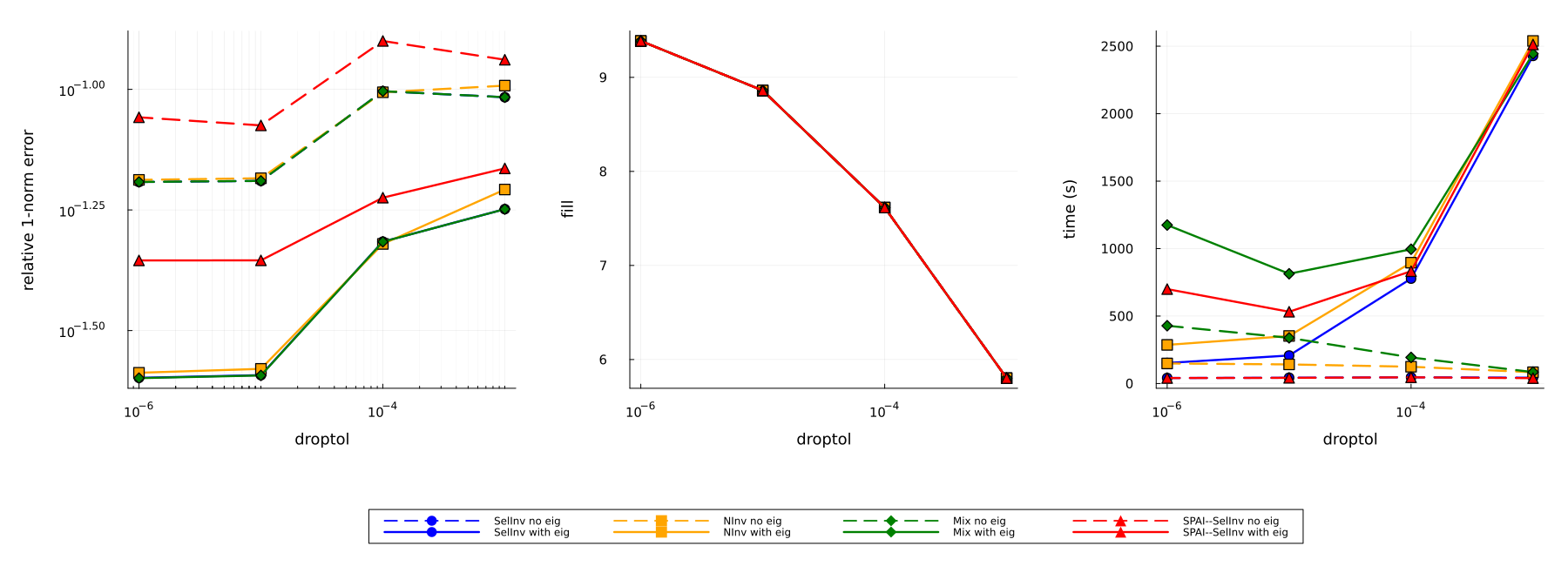}
    \caption{Detailed comparison of SelInv, NInv, Mix, and Mix--SPAI
    for the matrix \texttt{pwtk}. The panels show the relative
    diagonal error with and without spectral correction, fill from the ILDL factor, and
    the runtime components as functions of the drop tolerance.}
    \label{fig:pwtk_detailed}
\end{figure}

In the initial experiments, we compare the four base approximations $\AinvSelInv$, $\AinvNInv$, $\AinvMix$, and $A^{-1}_{\mathrm{Mix\mbox{-}SPAI}}$.
Our target here is to compare the four base approximations with and without the spectral correction. Mix-SPAI extends the Mix construction by incorporating an
additional sparse approximate-inverse refinement. The comparison is based on diagonal accuracy, fill, and runtime. The diagonal approximation generally depends on the drop tolerance $\tau$. As $\tau$ is tightened from $10^{-3}$ to $10^{-6}$, the ILDL factors generally become denser, and the base diagonal approximation becomes more accurate. The effect is particularly
focused in \texttt{crankseg\_2}, where the errors of SelInv, NInv, and Mix decrease strongly as $\tau$ is reduced, as seen in Figure \ref{fig:crankseg2_detailed}. The behavior of \texttt{pwtk}, shown in Figure~\ref{fig:pwtk_detailed}, is qualitatively similar, although the
relative benefit depends on the inverse construction and on the chosen
drop tolerance. In contrast, \texttt{oilpan} and \texttt{wathen120} show a much weaker dependence of the diagonal error on $\tau$; see the reference Table~\ref{tab:all_matrix_accuracy_cost}. For these matrices, the spectral correction is also largely ineffective, so the curves before and after correction nearly coincide.

The fill curves coincide for the four methods because they all use the same block ILDL factorization for a fixed value of $\tau$. Therefore, the factor fill is controlled by the drop tolerance and is independent of the subsequent inverse construction. Decreasing $\tau$ increases the factor fill because more entries are retained in the incomplete factors.

\begin{figure}[!t]
    \centering
    \includegraphics[width=\textwidth]{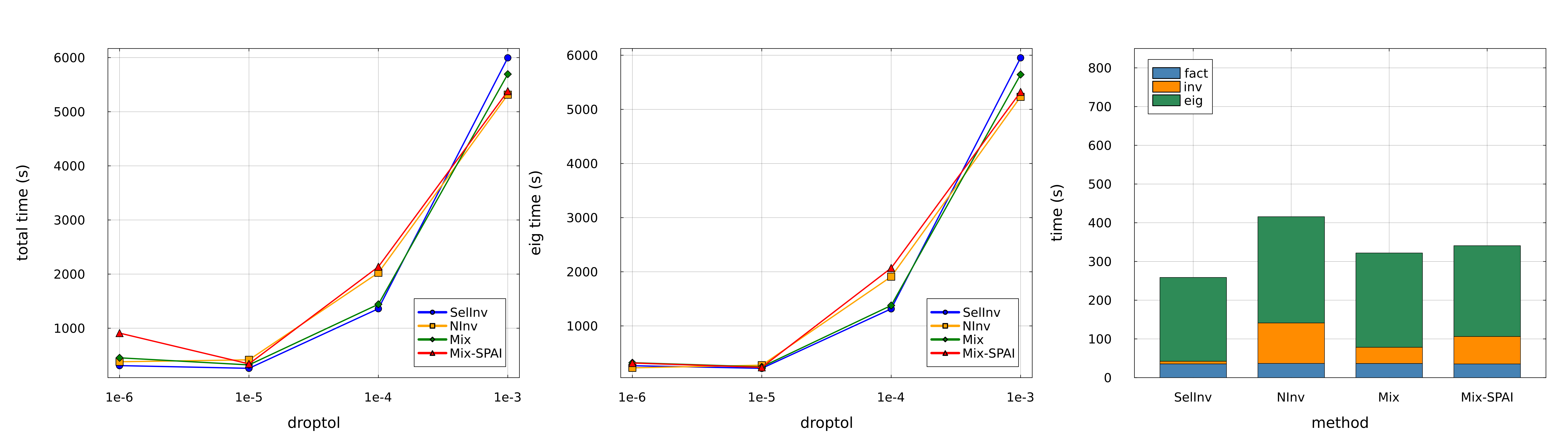}
   \caption{Runtime-focused comparison for the matrix \texttt{pwtk}. The left and middle panels show the total runtime and eigensolver runtime, respectively, as functions of the ILDL drop tolerance. The right panel decomposes the runtime into ILDL factorization, inverse construction, and eigensolver contributions at
$\tau=10^{-5}$.}
    \label{fig:pwtk_runtime_focus}
\end{figure}

The runtime does not depend on the ILDL factorization alone. It includes the common factorization cost, the method-dependent inverse construction cost, and the PRIMME eigensolver cost. The runtime results, as shown in Figures \ref{fig:pwtk_runtime_focus} and \ref{fig:oilpan_runtime_focus} illustrate that the spectral stage can constitute a substantial fraction of the total runtime even when no candidate correction is ultimately
accepted. Hence, the choice of $\tau$ should be based on an accuracy--cost tradeoff rather than on factor fill alone.

\subsection{Performance on Fixed Matrices}
\label{subsec: fixed matrix}
We evaluate all methods on the test set using SuiteSparse matrices, including structural engineering matrices such as \texttt{crankseg\_2}, \texttt{pwtk}, and \texttt{oilpan}. The underlying factorization framework is not restricted to real symmetric positive-definite (SPD) matrices and also supports symmetric indefinite and Hermitian problems. The characteristics of the corresponding matrix properties are summarized in Table \ref{tab:testmatrices}, while detailed
accuracy, fill, and runtime values at $\tau=10^{-5}$ are listed in Table \ref{tab:all_matrix_accuracy_cost}. The sizes of the matrices range from $19\,366$ to $217\,918$, as shown in Table \ref{tab:testmatrices}. From our chosen experiments, we observe that the necessity of spectral correction depends not only on the matrix but also on the choice of basis for the constructed inverse approximation. In particular, the same residual error may lead to different improvements: one for the inverse construction and another for accepting no improvement. We therefore present the spectral component as an optional correction whose benefit must be evaluated relative to the target quantity.

The corrected runtime decomposition also clarifies the relative cost of the four inverse constructions. SelInv requires only the selected inverse stage after the common ILDL factorization. NInv additionally constructs the factorized Neumann approximate inverse with the tighter inverse tolerance $\delta=10^{-2}$. The Mix method uses the cheaper NInv approximation with $\delta=10^{-1}$ together with the selected inverse. Mix--SPAI uses the same mix construction and subsequently applies SPAI refinement. The PRIMME times are of the same overall order for the different inverse approximations, although they need not be identical. Each application of the residual inverse-error operator

\[
\Eop=A^{-1}-\AinvCurrent.
\]

includes both an ILDL-preconditioned Krylov solve and an implementation of the method-dependent $A^{\mathrm{inv}}$. Therefore, the differences in sparsity and the one-time cost of applying $A^{\mathrm{inv}}$ depend on the executed total runtime of the eigensolver. In particular, using a tight drop tolerance in NInv may increase fill, thereby increasing the cost of a residual operator application.

\begin{figure}[!t]
    \centering
    \includegraphics[width=\textwidth]{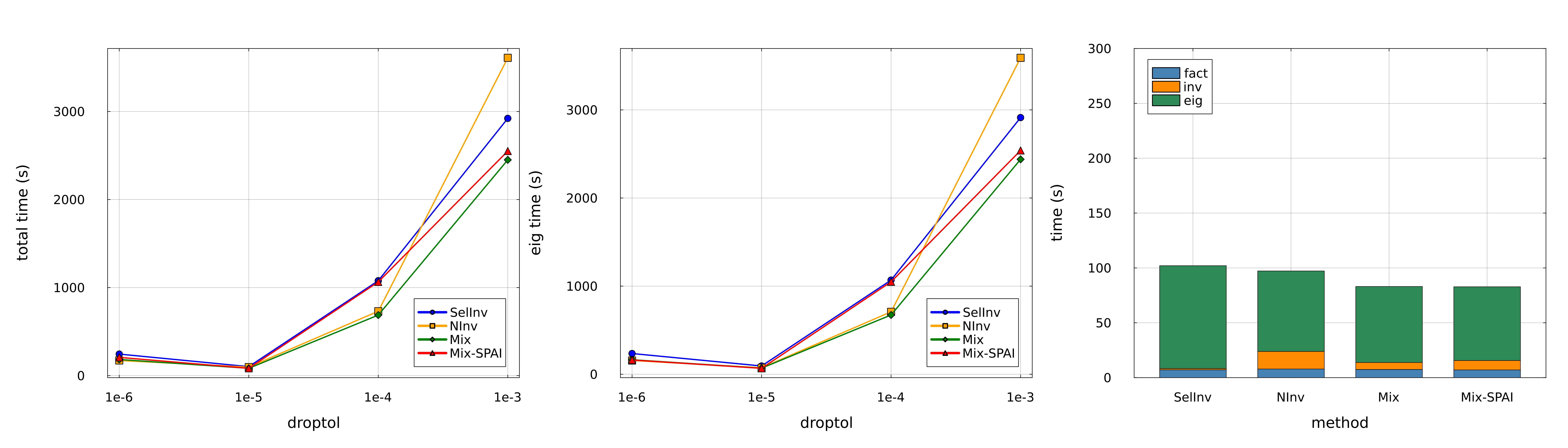}
   
     \caption{Execution time-focused comparison for the matrix \texttt{oilpan}. The left and middle plots show the total runtime and eigensolver runtime, respectively, as functions of the ILDL drop tolerance. The right panel decomposes the runtime into ILDL factorization, inverse construction, and eigensolver contributions at
$\tau=10^{-5}$.}
    \label{fig:oilpan_runtime_focus}
\end{figure}

The matrices \texttt{pwtk} and \texttt{apache1} represent a correction-helpful regime. For this larger problem, the spectral correction can substantially improve the diagonal approximation, although its effectiveness still depends on the underlying base inverse construction. For \texttt{pwtk} at $\tau=10^{-5}$, SelInv reduces the relative diagonal $1$-norm error from $6.46\times10^{-2}$ to $2.56\times10^{-2}$, corresponding to a gain of $2.53$,  with all three candidate eigenpairs accepted. NInv shows a comparable improvement, from $6.54\times10^{-2}$ to $2.63\times10^{-2}$, with a gain of $2.48$ and three accepted eigenpairs. By contrast, no spectral correction is accepted for Mix, so its final error remains $6.46\times10^{-2}$. Mix--SPAI accepts one correction mode and reduces the error from $3.37\times10^{-1}$ to $1.90\times10^{-1}$, corresponding to a gain of $1.77$. A similar behavior is observed for \texttt{apache1}. SelInv strengthens from $3.33\times10^{-1}$ to $1.69\times10^{-1}$,  while giving a gain of $1.97$ with the three accepted eigenpairs and NInv improves from $3.41\times10^{-1}$ to $1.88\times10^{-1}$ with a gain of $1.81$ and with the same number of accepted eigenpairs. The hybrid constructions also benefit moderately such that Mix reduces the error from $3.33\times10^{-1}$ to $2.27\times10^{-1}$ with one accepted mode, whereas Mix--SPAI decreases the error from $3.28\times10^{-1}$ to $2.30\times10^{-1}$, also with one accepted mode.

The matrix \texttt{crankseg\_2} illustrates a regime in which the computed spectral correction does not improve the target quantity. For SelInv at $\tau=10^{-5}$, one eigenpair satisfies the eigensolver accuracy and magnitude criteria, but the candidate diagonal error remains the same. Hence, the retained residual mode is not beneficial for the diagonal quantity of interest. For the benchmark matrices, where a reference diagonal is available, we apply a safeguard and retain the candidate only when
\[
\|\mathrm{error}_{\mathrm{cand}}\|_1
<
\|\mathrm{error}_{\mathrm{base}}\|_1.
\]
Otherwise, the candidate is rejected and the base approximation $\AinvCurrent$ is counted. Therefore, for these SelInv and Mix experiments, the final accepted error remains $5.07\times10^{-4}$. This result is theoretically important: a sufficiently accurate residual eigenpair does not, by itself, guarantee an improvement in a target-specific quantity. The detailed accuracy, fill, and runtime behavior is shown in Figure \ref{fig:crankseg2_detailed}.

For the matrices \texttt{oilpan} and \texttt{wathen120}, a no-change regime correction is demonstrated because of the negligibility or unacceptability of all computed eigenpairs. No correction modes are retained so that
\[
k_{\mathrm{acc}}=0
\]

 for each of the four inverse constructions. For \texttt{oilpan}, SelInv, NInv, and Mix give very similar relative diagonal errors, approximately $4.1\times10^{-3}$, whereas Mix--SPAI is less accurate. A similar behavior is observed for \texttt{wathen120}, where SelInv, NInv, and Mix yield relative errors close to $1.9\times10^{-2}$ and no spectral improvement is observed. Despite the lack of an accuracy gain, the spectral stage accounts for a substantial fraction of the computational cost. For SelInv, the eigensolver accounts for approximately $91.9\%$ of the total runtime
for \texttt{oilpan} and $78.1\%$ for \texttt{wathen120}. Thus, these two matrices illustrate that computing dominant residual modes can add considerable cost without improving the selected quantity of interest. The details of the cost and corresponding runtime decomposition are given in Table \ref{tab:all_matrix_accuracy_cost} and in Figure \ref{fig:oilpan_runtime_focus}.

\begin{figure}[t]
    \centering
    \includegraphics[width=\textwidth]{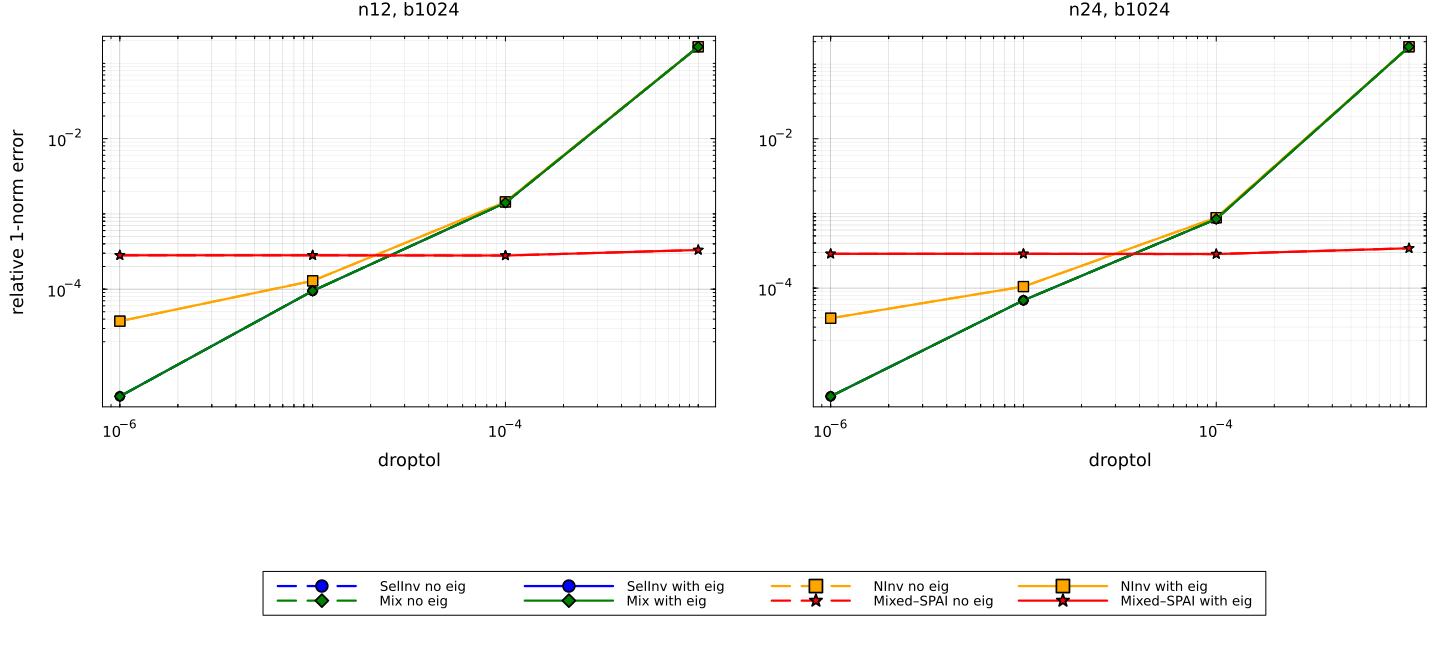}
    \caption{
    Relative 1-norm error of the approximate inverse diagonal for
    \texttt{NEW\_A\_microbench\_n12\_b1024} and
    \texttt{NEW\_A\_microbench\_n24\_b1024} as a function of the ILDL
    drop tolerance. Dashed curves represent the approximation without
    spectral correction, while solid curves include the PRIMME correction
    using three eigenpairs.
    }
    \label{fig:microbench-error-two-largest}
\end{figure}

The dense microbenchmark matrices confirm the same no-useful spectral-correction behavior observed for \texttt{oilpan} and \texttt{wathen120}. As shown in Table~\ref{tab:microbench-concise} and Figure~\ref{fig:microbench-error-two-largest}, SelInv and Mix achieve their smallest errors at $\tau=10^{-6}$, with relative diagonal errors of order $10^{-6}$, whereas NInv is about one order of magnitude less accurate but substantially sparser. Mix gives
approximately the sparsity of NInv while achieving SelInv-level accuracy. Mix--SPAI yields larger errors and additional computational cost, and no spectral correction is accepted in any reported case ($k_{\mathrm{acc}}=0$), indicating that the PRIMME stage provides no additional benefit for these matrices.

\begin{figure}[!t]
    \centering
    \includegraphics[width=\textwidth]{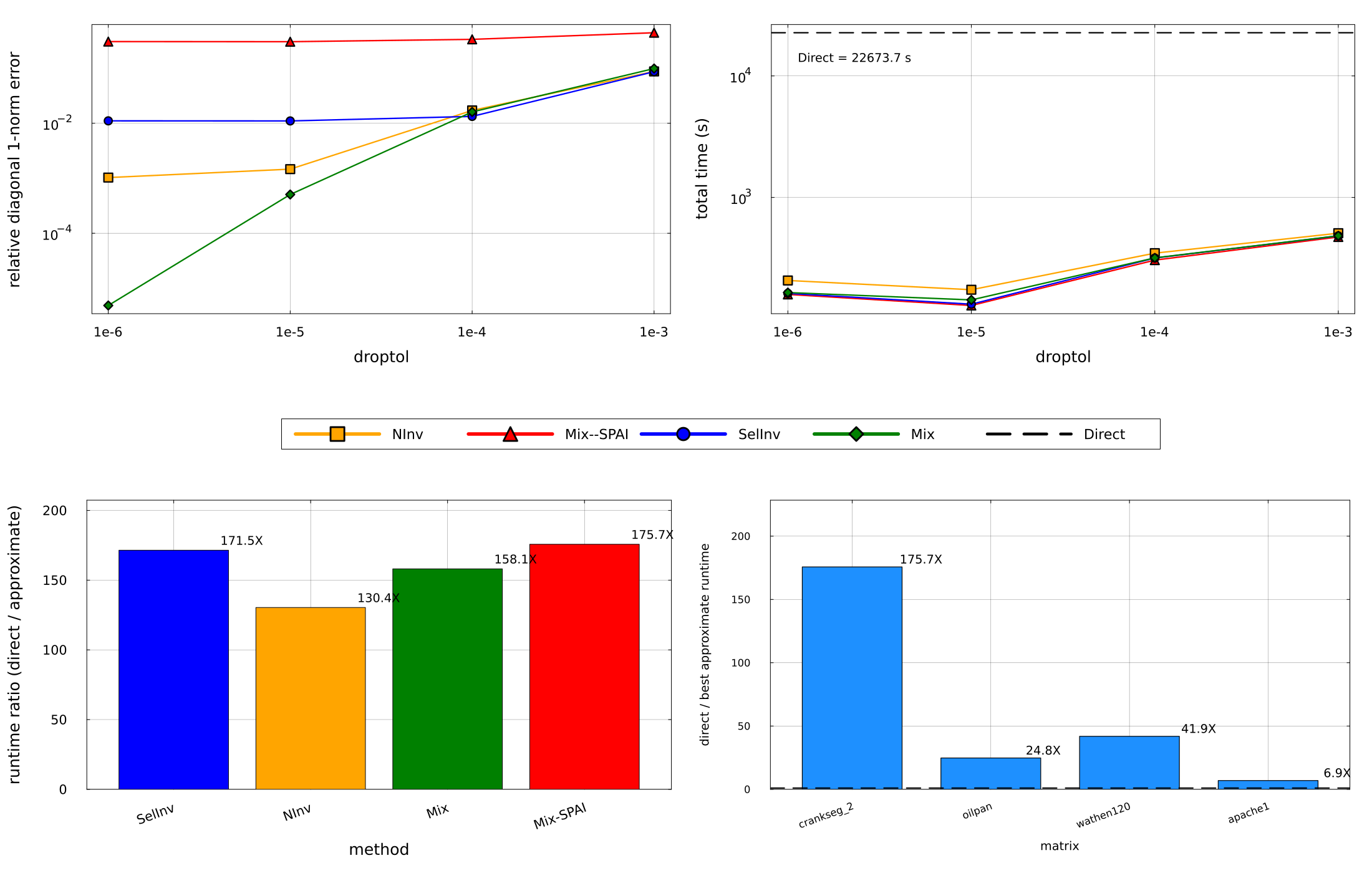}
    \caption{Comparison of the sparse direct reference with SelInv, NInv, Mix, and Mix--SPAI for the matrix \texttt{crankseg\_2}. The upper-left plot shows the relative diagonal error of the approximate methods versus the drop tolerance (with the direct result setting the reference error to zero). The upper-right plot shows the total approximate runtime together with the direct runtime. The lower-left plot compares the methods at $\tau=10^{-5}$, and the lower-right plot is a direct/best approximate runtime-focused comparison for the matrices: \texttt{crankseg\_2}, \texttt{oilpan}, \texttt{wathen120}, and \texttt{apache1}.}
    \label{fig:crankseg2_direct_vs_all_methods}
\end{figure}

\paragraph{Direct versus approximate inverse-diagonal computation.} The sparse direct method is considerably more expensive than the approximate methods. At $\tau=10^{-5}$, the fastest approximate
construction for \texttt{crankseg\_2} is about $176$ times faster than the sparse direct reference, as shown in Figure~\ref{fig:crankseg2_direct_vs_all_methods}. Across the four approximate methods, different accuracy--runtime trade-offs are observed, and a larger speedup does not necessarily refer a smaller
diagonal error.

The high cost of the direct computation does not arise only from the factorization, but also from the repeated triangular solves required to recover the complete inverse diagonal. Although the direct computation provides the reference solution with zero approximation error, it requires substantially greater computational effort. For the remaining matrices, the speedup of the approximate methods is also substantial, but its magnitude depends on both the matrix and the inverse
construction.

As the drop tolerance is tightened, the approximation error generally decreases, while the computational cost increases according to the inverse construction. Therefore, the methods should be assessed in
terms of the overall balance between diagonal accuracy and computational cost. These results demonstrate that approximate selected inversion can provide an efficient alternative to sparse direct
inverse-diagonal computation for large-scale problems.

\begin{table}[!t]
\centering
\small
\caption{Concise summary of the two largest successfully completed
microbenchmark experiments. For each matrix--method pair,
$\tau_{\mathrm{best}}$ denotes the ILDL drop tolerance that yields the
smallest final relative diagonal $1$-norm error over
$\tau\in\{10^{-6},10^{-5},10^{-4},10^{-3}\}$.
The inverse fill $f_{\mathrm{inv}}$ is measured relative to
$\operatorname{nnz}(A)$. The eigensolver share is defined as
$s_{\mathrm{eig}}=100T_{\mathrm{eig}}/T_{\mathrm{total}}$, and
$k_{\mathrm{acc}}$ denotes the number of accepted spectral-correction
modes.}
\label{tab:microbench-concise}

\begin{tabular}{llcccccc}
\toprule
Matrix & Method &
$\tau_{\mathrm{best}}$ &
$e_{\mathrm{final}}$ &
$f_{\mathrm{inv}}$ &
$T_{\mathrm{total}}$ [s] &
$s_{\mathrm{eig}}$ [\%] &
$k_{\mathrm{acc}}$ \\
\midrule

\multirow{4}{*}{\texttt{n12}}
& SelInv
& $10^{-6}$
& $3.76\times10^{-6}$
& $1.91$
& $96.17$
& $49.7$
& $0$ \\

& NInv
& $10^{-6}$
& $3.75\times10^{-5}$
& $5.09\times10^{-4}$
& $96.21$
& $63.5$
& $0$ \\

& Mix
& $10^{-6}$
& $3.76\times10^{-6}$
& $5.09\times10^{-4}$
& $104.82$
& $52.8$
& $0$ \\

& Mix--SPAI
& $10^{-4}$
& $2.80\times10^{-4}$
& $3.52\times10^{-3}$
& $105.70$
& $63.5$
& $0$ \\

\midrule

\multirow{4}{*}{\texttt{n24}}
& SelInv
& $10^{-6}$
& $3.54\times10^{-6}$
& $2.07$
& $259.88$
& $49.3$
& $0$ \\

& NInv
& $10^{-6}$
& $3.95\times10^{-5}$
& $4.98\times10^{-4}$
& $506.29$
& $72.1$
& $0$ \\

& Mix
& $10^{-6}$
& $3.54\times10^{-6}$
& $4.98\times10^{-4}$
& $531.32$
& $69.2$
& $0$ \\

& Mix--SPAI
& $10^{-4}$
& $2.85\times10^{-4}$
& $3.66\times10^{-3}$
& $579.93$
& $85.7$
& $0$ \\

\bottomrule
\end{tabular}

\vspace{1mm}
\begin{minipage}{0.96\textwidth}
\footnotesize
For SelInv,
$T_{\mathrm{inv}}=T_{\mathrm{SelInv}}$; for NInv,
$T_{\mathrm{inv}}=T_{\mathrm{NInv}}$; for Mix,
$T_{\mathrm{inv}}
=T_{\mathrm{SelInv}}+T_{\mathrm{NInv}}$; and for Mix--SPAI,
$T_{\mathrm{inv}}
=T_{\mathrm{SelInv}}+T_{\mathrm{NInv}}+T_{\mathrm{SPAI}}$.
For Mix--SPAI, SPAI refines the numerical values while keeping the
underlying Mix sparsity pattern fixed; thus, for a given Mix
approximation, the refinement itself does not introduce additional
nonzero entries. No spectral correction is accepted for any of the
reported microbenchmark cases, so that $e_{\mathrm{final}}=e_{\mathrm{base}}$ and
$k_{\mathrm{acc}}=0$ throughout. The different $f_{\mathrm{inv}}$ values for Mix and Mix--SPAI arise because their best results occur at different values of
$\tau_{\mathrm{best}}$.
\end{minipage}

\end{table}

\subsection{Application to covariance and uncertainty quantification.}

Let $Q\in\mathbb{R}^{n\times n}$ be a sparse precision matrix and let
\[
\Sigma=Q^{-1}
\]
be the corresponding covariance matrix. In many uncertainty-quantification applications, the full dense matrix $\Sigma$ is not required. Instead one needs selected quantities such as the diagonal variances $\Sigma_{ii}$, a small number of local covariances $\Sigma_{ij}$, or the trace
\[
\operatorname{tr}(\Sigma)=\sum_{i=1}^n\Sigma_{ii}.
\]
In the notation of this paper, we identify $A\equiv Q$, so that
$A^{-1}\equiv\Sigma$. Let
\[
S\subset\{1,\ldots,n\}^2
\]
be the prescribed set of entries, for example, the diagonal together with selected near-neighbor pairs. The selected covariance map is
\[
\Sigma_{\rm Sel}=P_S(\Sigma),
\qquad
[\Sigma_{\rm Sel}]_{ij}
=
\begin{cases}
\Sigma_{ij}, & (i,j)\in S,\\
0, & \text{otherwise}.
\end{cases}
\]
Using the approximation constructed in this paper, we use
\[
\widehat\Sigma
\approx
\AinvCorr
=
\AinvCurrent+V_k\Lambda_kV_k^T,
\]
or simply $\widehat\Sigma\approx\AinvCurrent$ if the spectral correction
is not accepted. The practically selected covariance approximation is then
\[
\Sigma_{\rm Sel}\approx P_S(\widehat\Sigma).
\]
In particular,
\[
\operatorname{diag}(\Sigma)\approx \operatorname{diag}(\widehat\Sigma),
\qquad
\operatorname{tr}(\Sigma)\approx
\mathbf{1}^T\operatorname{diag}(\widehat\Sigma).
\]
A probing-based approach may also be used when the main interest is the
diagonal of the inverse; see \cite{tang_probing_2012}.

\subsection{Covariance-type example without deflation}

As a structured numerical example, we consider the covariance-type matrix \cite{tang_probing_2012}
\[
C=\left(c_{ij}\right)_{i,j},
\qquad
c_{ij}
=
\begin{cases}
\left(1-\dfrac{d(i,j)}{\alpha}\right)^\beta,
& d(i,j)\le \alpha,\\[4pt]
0, & \text{otherwise},
\end{cases}
\]
where $d(i,j)$ denotes the distance between grid points $i$ and $j$, $\alpha>0$ is a cutoff radius, and $\beta>0$ controls the decay. This matrix provides a covariance-type benchmark for selected inverse entries and diagonal estimation.

\begin{figure}[!t]
    \centering
    \includegraphics[width=\textwidth]
    {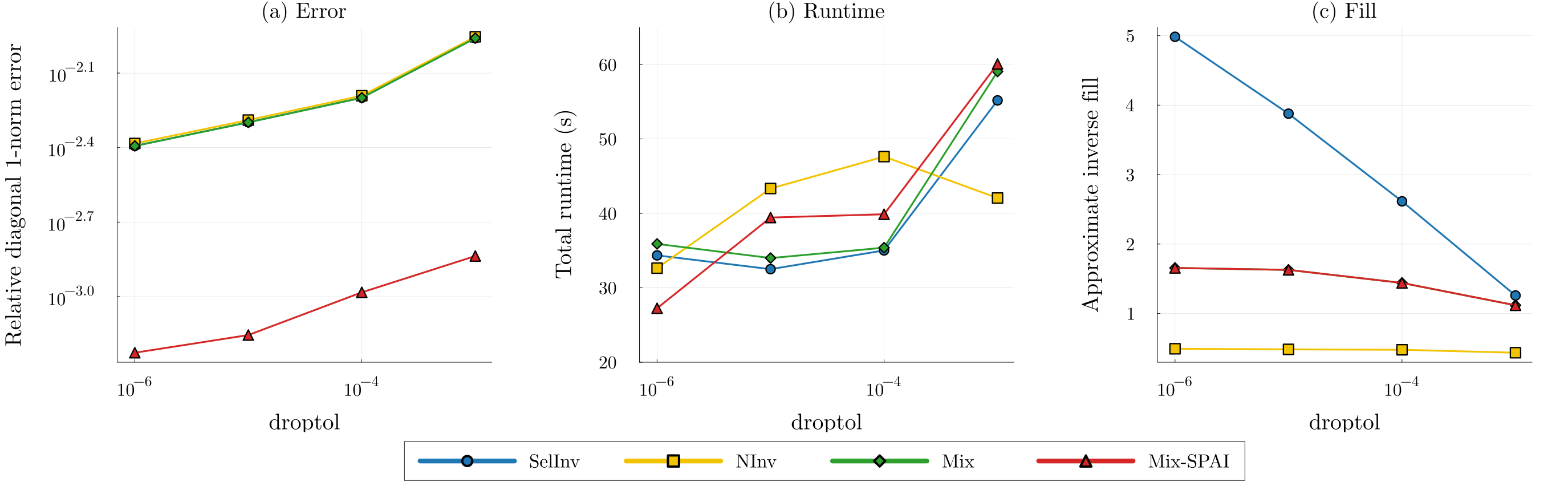}
    \caption{
    Detailed comparison of SelInv, NInv, Mix, and Mix-SPAI for the larger covariance matrix
    \texttt{covariance\_n200\_a3\_b4}, corresponding to matrix dimension $n=40000$.
    The panels show the relative diagonal $1$-norm error, total runtime, and approximate-inverse fill as functions of the drop tolerance.
    The total runtime includes the spectral check with three requested eigenpairs.
    No spectral correction was retained.
    }
    \label{fig:covariance_n200_a3_b4}
\end{figure}

\begin{figure}[!t]
    \centering
    \includegraphics[width=\textwidth]
    {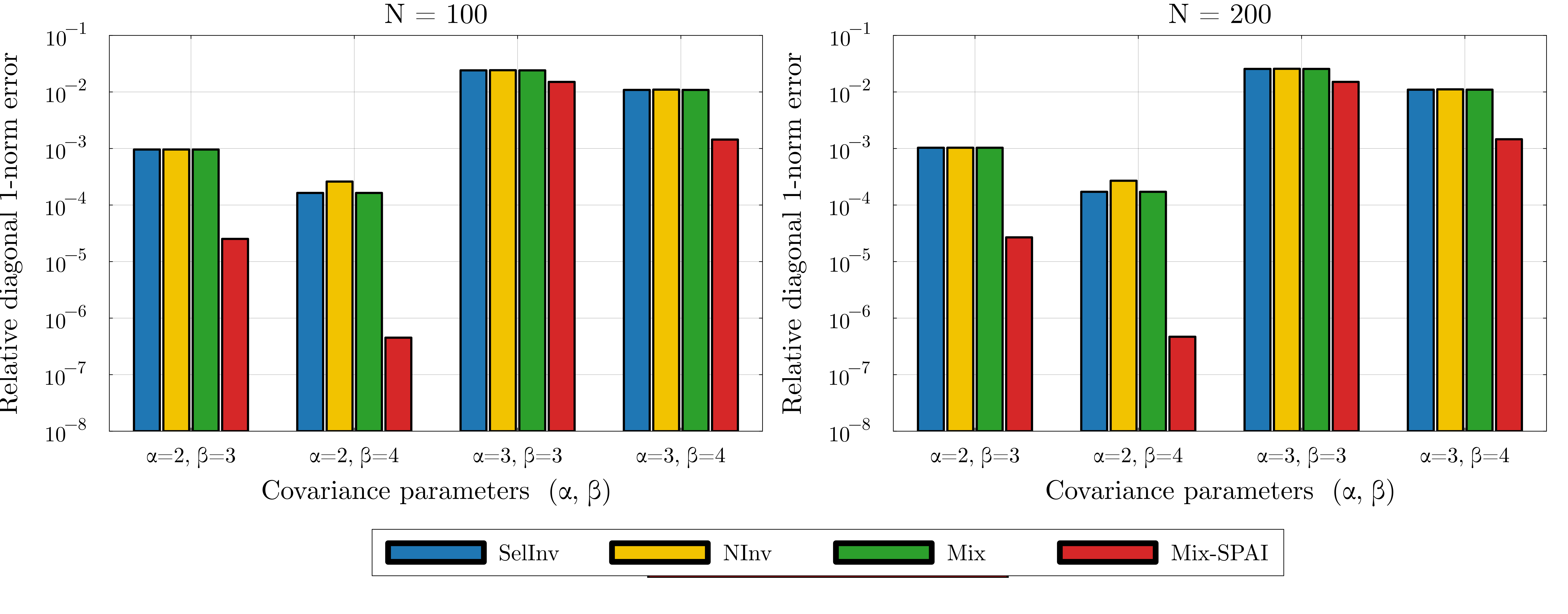}
    \caption{
    Relative diagonal $1$-norm errors for the covariance test
    problems at $\tau=10^{-3}$.
    Each $(N,\alpha,\beta)$ configuration defines a distinct
    covariance problem. For $N=100$ and $N=200$, the corresponding matrix dimensions
    are $n=10000$ and $n=40000$, respectively.
    Spectral correction was evaluated, but no correction was retained.
    Hence, the accepted errors coincide with the corresponding base errors.
    Mix-SPAI yields the smallest diagonal error for all displayed
    covariance cases, with particularly pronounced improvements
    for some parameter combinations.
    }
    \label{fig:covariance-error}
\end{figure}

In these covariance cases, the numerical behavior differs from that of the structural matrices above. The base approximations are already accurate, and the residual inverse-error spectrum showed no clear dominant low-rank effect. Spectral corrections with $3$, $5$, and $10$ requested eigenpairs were tested for the $N=200$ covariance problems, but no correction was retained. Hence, tightening the drop tolerance was more useful than applying spectral deflation.

Figure~\ref{fig:covariance_n200_a3_b4} illustrates this no-deflation scenario for the larger covariance matrix. The left panel shows the relative diagonal error, the middle panel shows the total runtime including the spectral check with three requested eigenpairs, and the right panel shows the approximate-inverse fill. For this larger covariance matrix, the spectral check adds substantial computational cost without a corresponding improvement in diagonal accuracy. This supports the adaptive strategy used in this work: the correction is applied only when its diagonal impact is significant.

To assess the impact of the covariance parameters, Figure~\ref{fig:covariance-error} compares the accepted diagonal errors for $N=100$ and $N=200$ at the fixed tolerance $\tau=10^{-3}$. Each $(N,\alpha,\beta)$ configuration defines a distinct covariance problem and is therefore shown as a separate category rather than connected by lines. Since no spectral correction is retained, the accepted errors coincide with the corresponding base errors. Mix-SPAI gives the smallest diagonal error for all displayed covariance cases, although the magnitude of the improvement depends on the covariance parameters.



\section{Concluding Remarks}
We have presented a scalable framework for the approximate selected inversion of large sparse symmetric systems, based on incomplete block $LDL^T$ factorizations. The numerical experiments show that the block ILDL drop tolerance is an important trade-off between accuracy, sparsity, runtime, and computational cost.

In this work, we considered four different inverse constructions, namely SelInv, NInv, Mixed, and Mix--SPAI, which provide different advantages depending on the matrix and the required accuracy. Spectral correction can further improve the approximation, but its effectiveness is highly matrix-dependent and must be balanced against the additional cost of the eigensolver. Therefore, spectral information is retained only when it yields a measurable improvement. For covariance-type matrices, the base approximation is already highly accurate, indicating that improving the underlying block  ILDL factorization can be more beneficial than applying a low-rank correction.

Overall, the proposed framework provides an adjustive balance among accuracy, sparsity, runtime, and computational cost, and the combination of block ILDL factorization, selected inversion, sparse approximate inverse techniques, and optional spectral correction enables the approximation to adapt to the matrix's numerical behavior.


\subsection{Future Work}
The present work mainly considers a single-level block ILDL framework. A natural next step is to investigate inverse-based multilevel factorizations, where the norm of the inverse triangular factor $L^{-1}$ is controlled level-by-level \cite{bollhofer_multilevel_2006}. Such a multilevel construction as implemented in ILUPACK\footnote{ILUPACK: \url{http://ilupack.tu-bs.de/}} may provide a better balance between robustness, fill, and approximation quality and will be considered in future work.

Future work will also involve computing-target-oriented work for spectral correction when the exact diagonal of $A^{-1}$ is unavailable. Possible approaches will focus on fully parallel implementations of the block ILDL factorization, adaptive strategies for automatically selecting the drop tolerances and correction rank, and extensions to larger-scale applications, including covariance estimation, uncertainty quantification, and electronic-structure calculations.

Another direction is a more systematic study of matrix structure, including 
Hermitian indefinite systems. In particular, it would be useful to determine which structural and spectral indicators predict whether a residual mode is aligned with a selected target quantity. We also believe that further study is needed into the mimicking behavior of Approximate Selected Inverse methods.

\section*{Acknowledgements}

This work was supported by Technische Universit\"at Braunschweig and the German Research Foundation (DFG), grant 470857344, in cooperation with SNSF grant 204817 at Universit\`a della Svizzera italiana, Switzerland, within the project ''Numerical Algorithms, Fundamental Structures, and Scalable Technologies for Highly Scalable Systems.''

\bibliographystyle{plain}
\bibliography{sample_zotero}
\appendix
\section{Notation}

Table~\ref{tab:notation} summarizes the main notation used throughout
the paper and the abbreviations used throughout the paper are ILDL (incomplete $LDL^T$ factorization), SelInv (selected inversion), NInv (Neumann-based approximate inverse), SPAI (sparse approximate inverse),
PCG (preconditioned conjugate gradient), and PRIMME (preconditioned iterative multimethod eigensolver).

\begin{table}[htbp]
\centering
\footnotesize
\setlength{\tabcolsep}{4pt}
\renewcommand{\arraystretch}{1.08}

\caption{Main notation used in this work.}
\label{tab:notation}

\begin{tabular}{p{0.22\textwidth} p{0.70\textwidth}}
\toprule
\textbf{Symbol} & \textbf{Description} \\
\midrule

$A$
& Sparse symmetric nonsingular system matrix. \\

$A^{-1}$
& Exact inverse of $A$. \\

$P$
& Symmetric permutation matrix used in the block ILDL factorization. \\

$L,D$
& Block unit lower-triangular and block-diagonal factors in
$A\approx LDL^T$. \\

$\Sset$
& Selected-entry pattern induced by the sparsity of $L+D+L^T$. \\

\midrule

$\AinvSelInv$
& Selected-inverse approximation on $\Sset$. \\

$\AinvNInv$
& Neumann-based sparse approximate inverse. \\

$\AinvMix$
& Mix approximation combining NInv and selected SelInv entries. \\

$A^{-1}_{\mathrm{Mix\text{-}SPAI}}$
& Mix approximation followed by SPAI refinement on the fixed pattern. \\

$\AinvCurrent$
& Generic base inverse approximation. \\

$\AinvCorr$
& Rank-$k$ spectrally corrected approximation of $A^{-1}$,
$\AinvCorr=\AinvCurrent+V_k\Lambda_kV_k^T$. \\
\midrule

$E_{\rm fact}$
& Factorization error, $E_{\rm fact}=A-LDL^T$. \\

$\Eop$
& Remaining inverse error,
$\Eop=A^{-1}-\AinvCurrent$. \\

$\tau$
& Drop tolerance of the block ILDL factorization. \\

$\delta$
& Inverse tolerance used in the Neumann-based approximation. \\

$\varepsilon$
& Threshold controlling the SelInv entries inserted into Mix. \\

\midrule

$V_k$
& Matrix of dominant eigenvectors of $\Eop$. \\

$\Lambda_k$
& Diagonal matrix of the corresponding eigenvalues. \\

$k$
& Number of candidate spectral correction modes. \\

$k_{\rm acc}$
& Number of accepted correction modes. \\

\midrule

$d,\widehat d$
& Exact and approximate inverse diagonals:
$d=\operatorname{diag}(A^{-1})$ and
$\widehat d=\operatorname{diag}(\AinvCurrent)$. \\

$\operatorname{err}_1$
& Relative diagonal $1$-norm error. \\

$\operatorname{fill}_{\rm fact}$
& ILDL factor fill relative to $\operatorname{nnz}(A)$. \\

$\operatorname{fill}_{\rm inv}$
& Inverse-approximation fill relative to $\operatorname{nnz}(A)$. \\

$T_{\rm fact},T_{\rm inv},T_{\rm eig}$
& Factorization, inverse-construction, and eigensolver times. \\

$T_{\rm total}$
& Total runtime,
$T_{\rm fact}+T_{\rm inv}+T_{\rm eig}$. \\

\bottomrule
\end{tabular}
\end{table}


      \begin{table*}[t]
\centering
\small
\setlength{\tabcolsep}{4.0pt}
\caption{Accuracy and computational-cost comparison at $\tau=10^{-5}$ and $n_{\mathrm{eig}}=3$. NInv uses inverse tolerance $\delta=10^{-2}$, whereas Mix and Mix--SPAI use the looser NInv tolerance $\delta=10^{-1}$. Here $e_{\mathrm{base}}$ is the relative diagonal $1$-norm error before spectral correction and $e_{\mathrm{final}}$ is the error after the benchmark acceptance safeguard. $G_{\mathrm{eig}}=e_{\mathrm{base}}/e_{\mathrm{final}}$, $k_{\mathrm{acc}}$ is the accepted correction rank, and $s_{\mathrm{eig}}=100T_{\mathrm{eig}}/T_{\mathrm{total}}$.}
\label{tab:all_matrix_accuracy_cost}
\resizebox{\textwidth}{!}{%
\begin{tabular}{llcccccccccc}
\toprule
Matrix & Method & $e_{\rm base}$ & $e_{\rm final}$ & $G_{\rm eig}$ & $k_{\rm acc}$ & $f_{\rm ILDL}$ & $T_{\rm ILDL}$ & $T_{\rm inv}$ & $T_{\rm eig}$ & $T_{\rm total}$ & $s_{\rm eig}$ (\%) \\
\midrule
crankseg\_2 & SelInv & $5.07\times10^{-4}$ & $5.07\times10^{-4}$ & 1.00 & 0 & 8.73 & 47.21 & 13.12 & 71.90 & 132.23 & 54.4 \\
crankseg\_2 & NInv & $1.47\times10^{-3}$ & $1.47\times10^{-3}$ & 1.00 & 0 & 8.73 & 48.47 & 48.71 & 76.74 & 173.93 & 44.1 \\
crankseg\_2 & Mix & $5.07\times10^{-4}$ & $5.07\times10^{-4}$ & 1.00 & 0 & 8.73 & 52.91 & 24.37 & 66.13 & 143.41 & 46.1 \\
crankseg\_2 & Mix--SPAI & $3.03\times10^{-1}$ & $3.03\times10^{-1}$ & 1.00 & 0 & 8.73 & 45.11 & 22.85 & 61.07 & 129.03 & 47.3 \\
\addlinespace[2pt]
pwtk & SelInv & $6.46\times10^{-2}$ & $2.56\times10^{-2}$ & 2.53 & 3 & 8.86 & 35.55 & 6.72 & 216.55 & 258.83 & 83.7 \\
pwtk & NInv & $6.54\times10^{-2}$ & $2.63\times10^{-2}$ & 2.48 & 3 & 8.86 & 36.73 & 104.53 & 274.38 & 415.65 & 66.0 \\
pwtk & Mix & $6.46\times10^{-2}$ & $6.46\times10^{-2}$ & 1.00 & 0 & 8.86 & 36.59 & 41.88 & 243.52 & 322.00 & 75.6 \\
pwtk & Mix--SPAI & $3.37\times10^{-1}$ & $1.90\times10^{-1}$ & 1.77 & 1 & 8.86 & 35.34 & 70.96 & 234.53 & 340.83 & 68.8 \\
\addlinespace[2pt]
oilpan & SelInv & $4.05\times10^{-3}$ & $4.05\times10^{-3}$ & 1.00 & 0 & 5.68 & 7.13 & 1.15 & 93.77 & 102.05 & 91.9 \\
oilpan & NInv & $4.07\times10^{-3}$ & $4.07\times10^{-3}$ & 1.00 & 0 & 5.68 & 7.85 & 16.06 & 73.32 & 97.23 & 75.4 \\
oilpan & Mix & $4.05\times10^{-3}$ & $4.05\times10^{-3}$ & 1.00 & 0 & 5.68 & 7.44 & 6.37 & 69.29 & 83.09 & 83.4 \\
oilpan & Mix--SPAI & $4.10\times10^{-3}$ & $4.10\times10^{-3}$ & 1.00 & 0 & 5.68 & 7.01 & 8.62 & 67.21 & 82.84 & 81.1 \\
\addlinespace[2pt]
wathen120 & SelInv & $1.88\times10^{-2}$ & $1.88\times10^{-2}$ & 1.00 & 0 & 4.92 & 2.17 & 0.23 & 8.56 & 10.96 & 78.1 \\
wathen120 & NInv & $1.89\times10^{-2}$ & $1.89\times10^{-2}$ & 1.00 & 0 & 4.92 & 2.23 & 0.53 & 5.57 & 8.33 & 66.9 \\
wathen120 & Mix & $1.88\times10^{-2}$ & $1.88\times10^{-2}$ & 1.00 & 0 & 4.92 & 2.18 & 0.43 & 8.61 & 11.22 & 76.7 \\
wathen120 & Mix--SPAI & $2.68\times10^{-1}$ & $2.68\times10^{-1}$ & 1.00 & 0 & 4.92 & 2.15 & 0.46 & 5.74 & 8.34 & 68.7 \\
\addlinespace[2pt]
apache1 & SelInv & $3.33\times10^{-1}$ & $1.69\times10^{-1}$ & 1.97 & 3 & 20.60 & 21.68 & 2.64 & 1125.03 & 1149.35 & 97.9 \\
apache1 & NInv & $3.41\times10^{-1}$ & $1.88\times10^{-1}$ & 1.81 & 3 & 20.60 & 20.85 & 63.60 & 1031.16 & 1115.62 & 92.4 \\
apache1 & Mix & $3.33\times10^{-1}$ & $2.27\times10^{-1}$ & 1.47 & 1 & 20.60 & 22.08 & 34.06 & 842.00 & 898.15 & 93.7 \\
apache1 & Mix--SPAI & $3.28\times10^{-1}$ & $2.30\times10^{-1}$ & 1.43 & 1 & 20.60 & 21.71 & 45.04 & 1111.88 & 1178.62 & 94.3 \\
\bottomrule
\end{tabular}%
}
\vspace{1mm}
\begin{minipage}{0.98\textwidth}
\footnotesize
For SelInv, $T_{\rm inv}=T_{\rm SelInv}$; for NInv, $T_{\rm inv}=T_{\rm NInv}$; for Mix, $T_{\rm inv}=T_{\rm SelInv}+T_{\rm NInv}$; and for Mix--SPAI, $T_{\rm inv}=T_{\rm SelInv}+T_{\rm NInv}+T_{\rm SPAI}$. The total runtime is reconstructed consistently as $T_{\rm total}=T_{\rm fact}+T_{\rm inv}+T_{\rm eig}$. A spectral candidate is accepted only when it decreases the relative diagonal $1$-norm error; otherwise $e_{\rm final}=e_{\rm base}$ and $k_{\rm acc}=0$.
\end{minipage}
\end{table*}
  


\end{document}